\documentclass{article}

\usepackage{standalone}
\usepackage[english]{babel}
\usepackage[utf8x]{inputenc}
\usepackage[numbers]{natbib}
\usepackage{authblk}
\usepackage[inline]{enumitem}

\usepackage{mathtools}
\usepackage{amssymb}
\usepackage{amsthm}
\usepackage{nicematrix}

\usepackage{tabularray}
\UseTblrLibrary{booktabs}

\usepackage[caption=false]{subfig}

\usepackage{algorithm}
\usepackage{algpseudocodex}

\algrenewcommand\algorithmicrequire{\textbf{Input:}}
\algrenewcommand\algorithmicensure{\textbf{Output:}}

\usepackage{tikz}
\usepackage{pgfplots}

\pgfplotsset{compat=1.18}
\usetikzlibrary{arrows.meta}
\usetikzlibrary{positioning}
\usetikzlibrary{calc}
\usetikzlibrary{math}
\usetikzlibrary{decorations}

\usepackage{hyperref}
\usepackage[noabbrev,nameinlink,capitalize]{cleveref}

\newtheorem{theorem}{Theorem}

\theoremstyle{definition}

\newtheorem{example}{Example}

\theoremstyle{remark}
\newtheorem{remark}[theorem]{Remark}

\hypersetup{
    colorlinks=true,
    linkcolor=blue,
    urlcolor=red,
    citecolor=teal,
}

\crefformat{prog}{#2Program~(#1)#3}
\Crefformat{prog}{#2Program~(#1)#3}
\crefmultiformat{prog}{#2Programs~(#1)#3}{ and~#2(#1)#3}{, #2(#1)#3}{, and~#2(#1)#3}
\Crefmultiformat{prog}{#2Programs~(#1)#3}{ and~#2(#1)#3}{, #2(#1)#3}{, and~#2(#1)#3}

\crefformat{cons}{#2Constraint~(#1)#3}
\Crefformat{cons}{#2Constraint~(#1)#3}
\crefmultiformat{cons}{#2Constraints~(#1)#3}{ and~#2(#1)#3}{, #2(#1)#3}{, and~#2(#1)#3}
\Crefmultiformat{cons}{#2Constraints~(#1)#3}{ and~#2(#1)#3}{, #2(#1)#3}{, and~#2(#1)#3}

\newcommand{\tr}[1]{#1^\top}

\newcommand{\midBig}{~\Big|~}

\newcommand{\setreal}{\mathbb{R}}
\newcommand{\setbin}{\{0, 1\}}

\newcommand{\vecone}{\mathbf{1}}

\newcommand{\setI}{\mathcal{I}}
\newcommand{\setIa}{\mathcal{I}_{1}}
\newcommand{\setIb}{\mathcal{I}_{2}}
\newcommand{\setT}{\mathcal{T}}
\newcommand{\setX}{\mathcal{X}}

\newcommand{\setV}{\mathcal{V}}
\newcommand{\setA}{\mathcal{A}}
\newcommand{\setB}{\mathcal{B}}

\newcommand{\setG}{\mathcal{G}}
\newcommand{\setE}{\mathcal{E}}

\newcommand{\bigO}{\mathcal{O}}

\newcommand{\funcR}{\mathbf{R}}

\newcommand{\vi}{v_i}
\newcommand{\ti}{t_i}
\newcommand{\tx}{\tr{t}x}
\newcommand{\vx}{\tr{v}x}
\newcommand{\vtx}{\tr{(v - t)}x}
\newcommand{\wx}{\tr{w}x}

\newcommand{\hx}{\hat{x}}
\newcommand{\zx}{z_{\hx}}
\newcommand{\vthx}{\tr{(v - t)}\hx}
\newcommand{\sumx}{\sum_{\hx\in\setX}}

\newcommand{\sk}{s_k}

\newcommand{\ba}{\bar{a}}
\newcommand{\ha}{\hat{a}}

\newcommand{\sumI}{\sum_{i\in\setI}}

\pgfdeclaredecoration{dashsoliddouble}{initial}{
  \state{initial}[width=\pgfdecoratedinputsegmentlength]{
    \pgfmathsetlengthmacro\lw{1.0pt+.5\pgflinewidth}
    \begin{pgfscope}
      \pgfpathmoveto{\pgfpoint{0pt}{\lw}}%
      \pgfpathlineto{\pgfpoint{\pgfdecoratedinputsegmentlength}{\lw}}%
      \pgfmathtruncatemacro\dashnum{%
        round((\pgfdecoratedinputsegmentlength-3pt)/6pt)
      }
      \pgfmathsetmacro\dashscale{%
        \pgfdecoratedinputsegmentlength/(\dashnum*6pt + 3pt)
      }
      \pgfmathsetlengthmacro\dashunit{3pt*\dashscale}
      \pgfsetdash{{\dashunit}{\dashunit}}{0pt}
      \pgfusepath{stroke}
      \pgfsetdash{}{0pt}
      \pgfpathmoveto{\pgfpoint{0pt}{-\lw}}%
      \pgfpathlineto{\pgfpoint{\pgfdecoratedinputsegmentlength}{-\lw}}%
      \pgfusepath{stroke}
    \end{pgfscope}
  }
}

\newcommand{\ie}{{\it i.e.~}}

\newcommand{\st}{\text{s.t.}}

\SetTblrTemplate{caption}{empty}

\title{Solving Combinatorial Pricing Problems using Embedded Dynamic Programming Models}
\author[1]{Quang Minh Bui}
\author[1]{Margarida Carvalho}
\author[2]{Jos\'e Neto}
\affil[1]{CIRRELT and D\'epartement d'informatique et de recherche op\'erationnelle, Universit\'e de Montr\'eal}
\affil[2]{T\'el\'ecom SudParis, Institut Polytechnique de Paris}
\date{}

\begin{document}

\maketitle

\begin{abstract}
  The combinatorial pricing problem (CPP) is a bilevel problem in which the leader
   maximizes their revenue by imposing tolls on certain items that they can control.
   Based on the tolls set by the leader, the follower selects a subset of items
   corresponding to an optimal solution of a combinatorial optimization problem.
   To accomplish the leader's goal, the tolls need to be sufficiently low to discourage
   the follower from choosing the items offered by the competitors.
   In this paper, we derive a single-level reformulation for the CPP by rewriting the follower's problem
   as a longest path problem using a dynamic programming model, and then taking its dual and applying strong duality.
   We proceed to solve the reformulation in a dynamic fashion with a cutting plane method.
   We apply this methodology to two distinct dynamic programming models, namely, a novel formulation designated as selection diagram and the well-known decision diagram. We also produce numerical results to evaluate their performances across three different specializations of the CPP
   and a closely related problem that is the knapsack interdiction problem. Our results showcase the potential of the two proposed reformulations over the natural value function approach, expanding the set of tools to solve combinatorial bilevel programs. 
\end{abstract}

\section{Introduction}
\label{sec:intro}
Consider a bilevel problem with two decision makers where one agent, the \textit{leader},
sets the \textit{tolls} (\ie markups) of certain items,
then the other agent, the \textit{follower}, selects a subset of items
according to another optimization problem (also known as the \textit{follower's problem}).
The follower has a choice between the items that the leader controls (called \textit{tolled items}),
and other items offered by the leader's competitors, whose prices, we assume, are fixed.
The goal of the leader is to maximize its revenue, which is the sum of tolls of the selected tolled items.
Although the leader naturally wants to set the tolls high,
such tolls cannot be excessive since the leader still needs the follower to choose their items in order to gain a profit.
Thus, the objectives of the two agents are neither adversarial nor cooperative.
In this paper, we investigate the case where the follower's problem is a combinatorial optimization problem
which can be expressed as a binary linear program.
We refer to the overall bilevel problem as \textit{combinatorial pricing problem} (CPP).
In the literature, the CPP is also called the Stackelberg pricing game.

\paragraph{Motivation}
The CPP is inspired by the network pricing problem (NPP) introduced by \citet{labbe1998}.
In that problem version, the NPP is a specific case of the CPP where the follower's problem is a shortest path problem.
In the CPP, the shortest path problem is replaced by other combinatorial problems such as
the minimum spanning tree problem \citep{cardinal2011}, the knapsack problem \citep{briest2012b, pferschy2021}, the stable set problem \citep{bohnlein2023},
and the bipartite vertex cover problem \citep{briest2012a}.
This class of problems is usually proven to be NP-hard \citep{cardinal2011,roch2005} or even $\Sigma_2^p$-hard \citep{bohnlein2023}.

Except for the NPP, most papers in the literature focus on the problems' computational complexity classification and approximation schemes.
Regarding the exact solution methods, general bilevel solvers like \texttt{MibS} \citep{tahernejad2020} or the one by \citet{fischetti2017} are not applicable due to
the leader's variables (the tolls in this case) being continuous and not appearing in the follower's constraints.
A popular method to solve the NPP is to find a dual of the follower's problem and then use the Karush-Kuhn-Tucker conditions
to convert the bilevel problem to a single-level reformulation \citep{bouhtou2007,bui2022,didibiha2006,heilporn2006}.
This method, referred to as \textit{dualize-and-combine}, is not trivial to generalize to the CPP because a dual of a binary linear program is not generally defined
and even if it exists, it is often computationally untractable.
A relatively novel technique to formulate a dual is described in \citet{lozano2022},
in which the authors use decision diagrams \citep{bergman2016} to derive
single-level reformulations for (bilevel) interdiction problems. This approach requires the assumption that the decision diagram model is of computationally tractable size allowing to apply off-the-shelf solvers directly to their reformulation.

\paragraph{Contributions}
Our first contribution is the application of the technique in \citet{lozano2022} to the CPP. 
We convert the follower's problem into an equivalent longest path problem corresponding to a dynamic programming model.
The longest path problem has a linear programming formulation,
hence we can dualize this linear program and combine it with the original primal follower's program to create a valid single-level reformulation for the CPP.
Besides decision diagrams, we also propose and evaluate a new dynamic programming model called \textit{selection diagram}.

Even with the reformulation well-defined, the computational untractability issue still remains due to the very likely large size of the reformulation. 
Thus, our second contribution is the proposition of a dynamic constraint generation scheme similar to the cutting plane method. Notably, we describe  methodologies to generate initial selection and decision diagrams, to add follower's solutions to them, and to avoid the increase on the number of layers of the decision diagrams.

Our third contribution is the numerical analysis comparing the performance of these formulations in the context of
three specializations of the CPP, in which the follower's problems are the knapsack problem, the maximum stable set problem,
and the minimum set cover problem, respectively.
We also apply the same methodology to the knapsack interdiction problem,
which can be reformulated as a special case of the CPP.
These results support a promising direction to scale the solving of the CPP.

Part of this work appears in the thesis of the first author~\citep{bui2024methods}.

\paragraph{Paper Organization}
In \cref{sec:comb-pricing}, we describe the bilevel formulation of the CPP and
its \textit{value function reformulation}, which serves as an introduction to the dualize-and-combine methodology.
We define two different dynamic programming models: selection diagram and decision diagram in \cref{sec:dp-model},
as well as the dynamic constraint generation algorithms corresponding to each model.
\Cref{sec:exp} presents the experimental setups and numerical results.
Finally, \cref{sec:conclusion} concludes the paper.

\section{Combinatorial Pricing Problem}
\label{sec:comb-pricing}

\subsection{Problem Description}
\label{ssec:desc}

Consider a combinatorial optimization problem with decision variables
$x\in\setX\subseteq\setbin^{\setI}$ indexed by the set of \textit{items} $\setI$.
Each item $i\in\setI$ has an associated \textit{base value} $\vi\in\setreal$.
The set $\setI$ is partitioned into the set of \textit{tolled items} $\setIa$
and the set of \textit{toll-free items} $\setIb$.
The leader wishes to impose a toll $\ti \geq 0$ for each tolled item $i\in\setIa$.
For the toll-free items $i\in\setIb$, we set $\ti = 0$.
Let $\setT = \setreal_{+}^{\setIa} \times \{0\}^{\setIb}$ be the set of feasible tolls.
The overall CPP is formulated as a bilevel program:
\begin{equation*}
	\max_{t,x} \{ \tx \mid t \in \setT,~x \in \funcR(t) \}
\end{equation*}
where $\funcR(t) \subseteq \setX$ is the set of optimal solutions of the
\textit{follower's problem} given the toll $t$, defined by:
\begin{equation}
	\label[prog]{prog:cpp-follower}
	\funcR(t) = \underset{x}{\arg\max} \{ \vtx \mid x \in \setX \}.
\end{equation}

Note that if $\funcR(t)$ is not a singleton, \ie there are multiple optimal solutions
of the follower's problem, then the follower will choose a solution $x$ that yields
the most revenue $\tx$ for the leader. This is commonly referred to as the
\textit{optimistic assumption}.

\begin{example}
	The \textit{knapsack pricing problem} (KPP) is a case of the CPP where
	the underlying optimization problem is a knapsack problem.
	Let $w \in \setreal_{+}^{\setI}$ and $C > 0$ be the knapsack weights and capacity,
	respectively. Then, the follower's problem of the KPP becomes:
	\begin{equation*}
		\funcR(t) = \underset{x}{\arg\max} \{ \vtx \mid \wx \leq C,~x \in \setbin^{\setI} \}.
	\end{equation*}
\end{example}

We make two remarks regarding the KPP:
\begin{remark}
	\label{rm:kpp-t-bound}
	The leader has no incentive to set $\ti > \vi$.
	Indeed, from the perspective of the follower's optimal value, $\vi - \ti < 0$ is the same as
	$\vi - \ti = 0$, and the former case does not give the leader any revenue,
	so the leader may as well set $\ti = \vi$.
	Thus, we can assume that $\ti \leq \vi$ for all $i \in \setI$.
\end{remark}

\begin{remark}
	\label{rm:kpp-maximal}
	There always exists an optimal solution $x^{*}$ of the KPP that is maximal.
	This follows from~\cref{rm:kpp-t-bound}, since including any new item
	produces a profit for both the follower and the leader. This property is a consequence of the \textit{monotonicity} of the lower-level feasible solutions~\citep{fischetti2019}.
\end{remark}

\subsection{Value Function Formulation}
\label{ssec:vf-form}

One method to solve the CPP is to convert the bilevel program into a single-level reformulation.
Such formulation may be obtained from a dual representation of the follower's problem \eqref{prog:cpp-follower}.
The follower's problems that we are interested in are binary programs,
hence the dual representation is not trivial to derive. Despite that, next we show how such
dual representation could be determined and used within a cutting plane method.

Define a new binary variable $\zx \in \setbin$ that indicates
if a particular solution $\hx \in \setX$ is selected or not.
We reformulate \cref{prog:cpp-follower} as follows:
\begin{equation}
	\label[prog]{prog:cpp-follower-z}
	\max_{z} \Big\{ \sumx (\vthx)\zx \midBig \sumx \zx = 1,~z \in \setbin^{\setX} \Big\}.
\end{equation}

\Cref{prog:cpp-follower-z} is totally unimodular,
hence we can relax $z \in \setbin^{\setX}$ to $z \geq 0$.
The dual of its linear relaxation is:
\begin{equation}
	\label[prog]{prog:cpp-follower-z-dual}
	\min_{L} \{ L \mid L \geq \vthx,~\hx \in \setX\}.
\end{equation}

To derive a single-level reformulation, in the bilevel program modeling the CPP, we replace the follower's optimality (\cref{prog:cpp-follower}) by its primal feasibility, dual feasibility (constraints of \cref{prog:cpp-follower-z-dual}) and an optimality condition, namely, the strong duality condition,
\begin{subequations}
	\label[prog]{prog:vf-single}
	\begin{align}
		\max_{t,x,L}\quad & \tx                                                                 \\
		\st\quad          & x \in \setX,                                                        \\
		                  & L \geq \vthx, &  &  \hx \in \setX, \label[cons]{cons:vf-single-dual} \\
		                  & \vtx = L,     &  & \label[cons]{cons:vf-single-sd}                  \\
		                  & t \in \setT.
	\end{align}
\end{subequations}
Observe that \cref{cons:vf-single-sd,cons:vf-single-dual} imply:
\begin{equation}
	\label{eq:vf-constraint}
	\vtx \geq \vthx, \qquad \hx \in \setX.
\end{equation}

This is usually referred as the \textit{(optimal) value function constraint}.
Thus, we refer to \cref{prog:vf-single} as the \textit{value function formulation}.

\begin{example}
	\label{ex:kpp-vf}
	Consider a KPP instance of four items $\setI = \{1,2,3,4\}$.
	The first three items are tolled, \ie $\setIa = \{1,2,3\}$, and $\setIb = \{4\}$.
	The base values are $v = (1,1,1,1)$. The knapsack weights and capacity are
	$w = (1,1,1,2)$ and $C = 3$ respectively. Thus, we have
	\[ \setX = \{x \in \setbin^{4} \mid x_{1} + x_{2} + x_{3} + 2x_{4} \leq 3\}. \]

	By \cref{rm:kpp-maximal}, we only need to consider maximal solutions for~\cref{cons:vf-single-dual}. These could be obtained by enumerating all feasible solutions and keeping the inclusionwise maximal. 
	There are four maximal solutions in $\setX$: $(1,1,1,0)$, $(1,0,0,1)$, $(0,1,0,1)$, and $(0,0,1,1)$.
	It is more convenient to refer to them as the sets they represent: $\{1,2,3\}$, $\{1,4\}$,
	$\{2,4\}$, and $\{3,4\}$. The value function formulation
	corresponding to this KPP instance is:
	\begin{align*}
		\max_{t,x,L}\quad & t_{1}x_{1} + t_{2}x_{2} + t_{3}x_{3}                                \\
		\st\quad          & x_{1} + x_{2} + x_{3} + 2x_{4} \leq 3,                              \\
		                  & L \geq 3 - t_{1} - t_{2} - t_{3},                                   \\
		                  & L \geq 2 - t_{1},                                                   \\
		                  & L \geq 2 - t_{2},                                                   \\
		                  & L \geq 2 - t_{3},                                                   \\
		                  & (1 - t_{1})x_{1} + (1 - t_{2})x_{2} + (1 - t_{3})x_{3} + x_{4} = L, \\
		                  & x \in \setbin^{4},                                                  \\
		                  & t \geq 0.
	\end{align*}
\end{example}

To solve \cref{prog:vf-single} as a mixed-integer linear program (MILP), we need to linearize the bilinear terms $\tx$.
This can be done using the McCormick envelope \citep{mccormick1976} by replacing $t_{i}x_{i}$ for each $i \in \setIa$
with $s_{i}$, accompanied by the following constraints:
\begin{align}\label{eq:mccormick}
	0 & \leq s_{i} \leq M_{i}x_{i},
	 & 0 & \leq t_{i} - s_{i} \leq M_{i}(1 - x_{i}),
\end{align}
where $M_{i}$ is a sufficiently large bound of $t_{i}x_{i}$.
For instance, a valid bound for the KPP is $M_{i} = \vi$ due to \cref{rm:kpp-t-bound}.

Typically, the size of $\setX$ is exponential, so in practice, we generate
\cref{cons:vf-single-dual} dynamically in a cutting plane fashion as described in \cref{alg:vf-dynamic}. In this algorithm, the (relaxed) master program corresponds to \cref{prog:vf-single} with a subset of Constraints~\eqref{cons:vf-single-dual} solved in line 3, and the subproblem corresponds to the follower's problem solved in line 4. The remaining steps correspond to verifying whether the master program provides an optimal solution to the follower. If yes, the method terminates. Otherwise, a new Constraint~\eqref{cons:vf-single-dual} is added to the master program.

\begin{algorithm}[ht]
	\caption{Solution of \cref{prog:vf-single} with dynamically generated constraints}
	\label{alg:vf-dynamic}

	\begin{algorithmic}[1]
		\Require{The CPP instance.}
		\Ensure{An optimal solution $(t^*, x^*)$ of the CPP instance.}
		\State Create the master problem from \cref{prog:vf-single} without \cref{cons:vf-single-dual}.
		\Loop
		\State Solve the master problem and let $(\tilde{t}, \tilde{x})$ be the current optimal solution.
		\State Solve \cref{prog:cpp-follower} with $t = \tilde{t}$ to get $\hx$.
		\Comment{\ie choose any $\hx \in \funcR(\tilde{t})$}

		\If{$\tr{(v-\tilde{t})}\tilde{x} < \tr{(v-\tilde{t})}\hat{x}$}
		\LComment{$(\tilde{t},\tilde{x})$ is not bilevel feasible. }
		\State Add the constraint $L \geq \vthx$ to the master program.

		\Else
		\LComment{$(\tilde{t},\tilde{x})$ is bilevel feasible. Stop.}
		\State $(t^*, x^*) \gets (\tilde{t}, \tilde{x})$
		\State \textbf{break}
		\EndIf
		\EndLoop
		\State \Return $(t^*, x^*)$
	\end{algorithmic}
\end{algorithm}

\section{Embedded Dynamic Programming Model}
\label{sec:dp-model}

The value function formulation in \cref{prog:vf-single} is simple and general.
However, the structure of the follower's problem is oblivious to the dual representation.
We aim to introduce some structural information back to the dual in the form of a dynamic programming (DP) model.
It is well-known that any DP model with a finite number of states and actions can be represented as a directed multigraph \citep{bertsekas2012}.
Thus, for our convenience, we define all DP models by their graph representations.
Suppose that the follower's problem possesses a DP model whose definition includes:
\begin{itemize}
	\item A set of \textit{nodes} $\setV$ (corresponding to the \textit{states});
	\item A set of \textit{arcs} $\setA$ (corresponding to the \textit{actions});
	\item Two functions $u, w : \setA \to \setV$ that return
	      the source node $u(a)$ and the target node $w(a)$ of arc $a$, respectively;
	\item An \textit{item function} $I: \setA \to 2^\setI$ associating to each arc $a\in \setA$ a subset of items $I(a)$, corresponding to items selected by the follower; 
	\item An \textit{initial node} $p \in \setV$ and a \textit{terminal node} $q \in \setV$.
\end{itemize}

The directed multigraph $(\setV, \setA, u, w)$ must be acyclic.
We define the \textit{length} $F(a; t)$ of arc $a \in \setA$ given toll $t \in \setT$ as
the sum of the tolled values of all the items in $I(a)$:
\begin{equation*}
	F(a; t) = \sum_{i \in I(a)} (v_i - t_i).
\end{equation*}

The objective of the DP model is to find the longest path from $p$ to $q$.
Henceforth, when we refer to \textit{paths}, we mean paths from $p$ to $q$.
This optimization problem can be expressed by the following formulation
(which is the dual of a typical linear formulation for the longest path problem):
\begin{subequations}
	\label[prog]{prog:dp-dual}
	\begin{align}
		\min_{y}\quad & y_p                                                 \\
		\st\quad      & y_{u(a)} \geq F(a; t) + y_{w(a)}, &  &  a \in \setA, \\
		              & y_q = 0.
	\end{align}
\end{subequations}

Because the DP model is equivalent to the follower's problem (\ie they return the same optimal objective value),
we can use \cref{prog:dp-dual} as a dual representation in a similar way to \cref{prog:cpp-follower-z-dual}.
The result is a new single-level reformulation\footnotemark:
\begin{subequations}
	\label[prog]{prog:dp-single}
	\begin{align}
		\max_{t,x,y}\quad & \tx                                                                                  \\
		\st\quad          & x \in \setX,     \label[cons]{cons:primalproblem}                                                                     \\
		                  & y_{u(a)} \geq F(a; t) + y_{w(a)}, &  &  a \in \setA, \label[cons]{cons:dp-single-arc} \\
		                  & y_q = 0,                          &  & \label[cons]{cons:dp-single-term}             \\
		                  & \vtx = y_p,                       &  &                                               \\
		                  & t \in \setT.
	\end{align}
\end{subequations}
\footnotetext{\cite{lozano2022}  obtain a similar single-level reformulation. However,
\cref{prog:dp-single} uses the follower's original representation of the feasible region, \cref{cons:primalproblem}, while
~\cite{lozano2022} uses the path formulation resulting from the DP model.}

\begin{figure}[ht]
	\centering
	\scalebox{0.9}{\begin{tikzpicture}
	\tikzstyle{state}=[draw,circle,inner sep=0pt,minimum size=24pt]
	\tikzstyle{action}=[circle,fill=white,inner sep=1pt]

	\node[state] (S1) at (0, 0) { $p$ };
    \node[state] (S2) at (3, 0) { $q$ };

    \draw (S1) to[bend left=80] node[action] {$\hx^1$} (S2);
    \draw (S1) to[bend left=40] node[action] {$\hx^2$} (S2);
    \draw (S1) --node[action] {$\hx^3$} (S2);
    \draw (S1) to[bend right=40] node[action] {$\hx^4$} (S2);
    \draw (S1) to[bend right=80] node[action] {$\hx^5$} (S2);

\end{tikzpicture}}
	\caption{Value function formulation as a DP model.}
	\label{fig:vfgraph}
\end{figure}

\begin{example} We show that the value function formulation \cref{prog:vf-single} can be thought
	as a special case of \cref{prog:dp-single}.
	Consider a DP model consisting of only two nodes $p$ and $q$.
	Each solution $\hx \in \setX$ corresponds to an arc $a_{\hx} \in \setA$.
	The item function is defined to be the set of items selected in $\hx$:
	\begin{equation*}
		I(a_{\hx}) = \{ i \in \setI \mid \hx_i = 1 \}.
	\end{equation*}
	Then, the variable $L$ in \cref{prog:vf-single} is identical to $y_p$ in \cref{prog:dp-single}. Thus, under this specific DP model, these programs are identical.
	The DP graph of the value function formulation is shown in \cref{fig:vfgraph}.
	In all graphs illustrated in this paper, the arcs are always directed from left to right. \cref{prog:dp-single} can also encompass other single-level formulations under different DP models, as we will see in the next sections.
\end{example}

Next, we explore the application of \cref{prog:dp-single} to two types of DP models:
\textit{selection diagrams} (\cref{ssec:sdiag}) and \textit{decision diagrams} (\cref{ssec:ddiag}).
The former defines the actions based on the items to be selected,
while the latter defines the actions as decisions to include or exclude a specific item.
For each type of model, we also describe the procedure to solve \cref{prog:dp-single}
dynamically as in \cref{alg:vf-dynamic} alongside any other technical improvement.

\subsection{Selection Diagram}
\label{ssec:sdiag}

The first family of dynamic programming models that we investigate is called \textit{selection diagrams}.
In a selection diagram, at each step, the follower either selects an item to be included into the current partial solution,
or decides to stop and returns the solution.
In particular, let $\setI(\hx) = \{i \in \setI \mid \hx_i = 1\}$ be the set of items selected in $\hx \in \setX$.
The set of nodes in a selection diagram is the family of subsets of $\setI$ that are contained in some feasible solution in $\setX$,
adjoined with a special terminal node $q$:
\begin{equation*}
	\setV = \{J \subseteq \setI \mid J \subseteq \setI(\hx) ~\text{for some}~ \hx \in \setX \} \cup \{q\}.
\end{equation*}
The initial node $p$ in this case is the node $\varnothing$.

For every pair of nonterminal nodes $J$, $K \in \setV \setminus \{q\}$ such that $J \subseteq K$ and $|K| - |J| = 1$,
we add an arc $a_{JK}$ from $J$ to $K$ such that $I(a_{JK}) = K \setminus J$.
Evidently, $I(a_{JK})$ is a singleton.
Moreover, given a nonterminal node $J$, if $J = \setI(\hx)$ for some $\hx \in \setX$,
then we connect $J$ directly to $q$ with a \textit{terminal arc} $\ba_J$ with $I(\ba_J) = \varnothing$.

With this construction, the nodes on a path from $\varnothing$ to $q$
form a sequence of nested subsets of $\setI$:
\begin{equation*}
	\varnothing ~,~ \{i_1\} ~,~ \{i_1, i_2\} ~,~ \ldots ~,~ \{i_1,\ldots,i_n\} ~,~ q
\end{equation*}
such that $\{i_1,\ldots,i_n\} = \setI(\hx)$ for some $\hx \in \setX$.
The length of such a path is exactly $\vthx$.
Thus, finding the longest path in a selection diagram is equivalent to finding the solution $\hx \in \setX$
that yields the highest objective value for the follower.

\begin{example}
	\label{ex:kpp-sd}
	The selection diagram for the KPP instance in \cref{ex:kpp-vf} is illustrated in \cref{fig:kpp-sd}. Recall that the only toll-free item is item 4.
	Once again, we only keep track of the four maximal solutions: $\{1,2,3\}$, $\{1,4\}$, $\{2,4\}$, $\{3,4\}$.
	Note that a given solution $\hx \in \setX$ may be represented by multiple paths from $\varnothing$ to $q$.
	Applying \cref{prog:dp-single} to \cref{fig:kpp-sd},
	\cref{cons:dp-single-arc,cons:dp-single-term} become:
	\begin{small}
		\begin{align*}
			y_{\varnothing} & \geq 1 - t_1 + y_1,     & y_{\varnothing} & \geq 1 - t_2 + y_2,     & y_{\varnothing} & \geq 1 - t_3 + y_3,     & y_{\varnothing} & \geq 1 + y_4,          \\
			y_1             & \geq 1 - t_2 + y_{12},  & y_2             & \geq 1 - t_1 + y_{12},  & y_3             & \geq 1 - t_1 + y_{13},  & y_4             & \geq 1 - t_1 + y_{14}, \\
			y_1             & \geq 1 - t_3 + y_{13},  & y_2             & \geq 1 - t_3 + y_{23},  & y_3             & \geq 1 - t_2 + y_{23},  & y_4             & \geq 1 - t_2 + y_{24}, \\
			y_1             & \geq 1 + y_{14},        & y_2             & \geq 1 + y_{24},        & y_3             & \geq 1 + y_{34},        & y_4             & \geq 1 - t_3 + y_{34}, \\
			y_{12}          & \geq 1 - t_3 + y_{123}, & y_{13}          & \geq 1 - t_2 + y_{123}, & y_{23}          & \geq 1 - t_1 + y_{123},                                            \\
			y_{14}          & \geq 0,                 & y_{24}          & \geq 0,                 & y_{34}          & \geq 0,                 & y_{123}         & \geq 0.
		\end{align*}
	\end{small}

	\begin{figure}[ht]
		\centering
		\scalebox{0.9}{\begin{tikzpicture}[xscale=2]
	\tikzstyle{state}=[draw,circle,inner sep=0pt,minimum size=24pt]

	\node[state] (S0) at (0, -1.5) { $\varnothing$ };

	\node[state] (S1) at (1, -0) { $\{1\}$ };
	\node[state] (S2) at (1, -1) { $\{2\}$ };
	\node[state] (S3) at (1, -2) { $\{3\}$ };
	\node[state] (S4) at (1, -3) { $\{4\}$ };

	\node[state] (S12) at (2, +1) { \small $\{1,2\}$ };
	\node[state] (S13) at (2, -0) { \small $\{1,3\}$ };
	\node[state] (S23) at (2, -1) { \small $\{2,3\}$ };
	\node[state] (S14) at (2, -2) { \small $\{1,4\}$ };
	\node[state] (S24) at (2, -3) { \small $\{2,4\}$ };
	\node[state] (S34) at (2, -4) { \small $\{3,4\}$ };

	\node[state] (S123) at (3, 0) { \small $\{1,2,3\}$ };

	\node[state] (BS) at (4, -1.5) { $q$ };

    \draw (S0) -- (S1);
    \draw (S0) -- (S2);
    \draw (S0) -- (S3);
    \draw (S0) -- (S4);

    \draw (S1) -- (S12);
    \draw (S1) -- (S13);
    \draw (S1) -- (S14);

    \draw (S2) -- (S12);
    \draw (S2) -- (S23);
    \draw (S2) -- (S24);

    \draw (S3) -- (S13);
    \draw (S3) -- (S23);
    \draw (S3) -- (S34);

    \draw (S4) -- (S14);
    \draw (S4) -- (S24);
    \draw (S4) -- (S34);

    \draw (S12) -- (S123);
    \draw (S23) -- (S123);
    \draw (S13) -- (S123);

    \draw[dashed] (S123) -- (BS);
    \draw[dashed] (S14) -- (BS);
    \draw[dashed] (S24) -- (BS);
    \draw[dashed] (S34) -- (BS);

\end{tikzpicture}}
		\caption{Selection diagram of \cref{ex:kpp-sd}.
			Dashed arcs represent the terminal arcs. }
		\label{fig:kpp-sd}
	\end{figure}
\end{example}

\subsubsection{Dynamic Generation}
\label{sssec:sdiag-dynamic}
The number of subsets of $\setI$ is exponential with respect to $|\setI|$, so in practice,
enumerating the full set of \cref{cons:dp-single-arc,cons:dp-single-term}, as in \cref{ex:kpp-sd}, is untractable.
Thus, a method to generate the variables and constraints dynamically is required.
Furthermore, some popular MILP solvers do not support column generation,
so the dynamic generation of the variables $y_J$ is usually not a viable option. 
In this section, we describe an algorithm to solve \cref{prog:dp-single}
using a selection diagram which is updated dynamically only with new constraints
(which are often referred to as \textit{lazy constraints}).

The first step is to choose a set of nodes $\hat\setV$ that will be used throughout the procedure.
This set of nodes is final, since we do not opt for column generation.
Evidently, the initial node $\varnothing$ and the terminal node $q$ must be in $\hat\setV$.
For the remaining nodes, we limit our options to only singletons and pairs,
\ie we only consider the nodes $J$ where $|J| \in \{1, 2\}$.
Even so, some pairs may not belong to a feasible solution, and the number of all pairs is often too large,
hence not every pair should be present in the selection diagram.
We propose a sampling algorithm described in \cref{alg:pair-sampling}.
The algorithm iteratively samples a solution $\hx \in \setX$. The solution sampling is problem-specific, and for problems satisfying the monotonicity property (\cref{rm:kpp-maximal}), sampled solutions must be maximal. Then, the algorithm randomly selects a pair within $\setI(\hx)$.
After reaching the desired number of pairs, it builds the connections within the first 2 layers.
The end result is an initial selection diagram to be used in \cref{prog:dp-single}.

\begin{algorithm}[ht]
	\caption{Generation of the initial selection diagram}
	\label{alg:pair-sampling}

	\begin{algorithmic}[1]
		\Require{The CPP instance, the desired number of pairs \(N\).}
		\Ensure{The set of nodes $\hat\setV$ and the set of initial arcs $\hat\setA$.}
		\LComment{Initial nodes}
		\State $\hat\setV_0 \gets \{\varnothing\}$, $\hat\setV_1 \gets \varnothing$, $\hat\setV_2 \gets \varnothing$
		\For{$n = 1,\ldots,N$}
		\State Sample $\hx \in \setX$
		\State Sample a random pair $\{i, j\}$ within $\setI(\hx)$
		\State $\hat\setV_1 \gets \hat\setV_1 \cup \{i, j\}$
		\Comment{Set of singletons}
		\State $\hat\setV_2 \gets \hat\setV_2 \cup \{\{i, j\}\}$
		\Comment{Set of pairs}
		\EndFor
		\State $\hat\setV \gets \hat\setV_0 \cup \hat\setV_1 \cup\hat\setV_2 \cup \{q\}$

		\LComment{Initial arcs}
		\State $\hat\setA \gets \varnothing$
		\For{$k = 1,2$}
		\ForAll{$J \in \hat\setV_{k-1}$, $K \in \hat\setV_k$}
		\If{$J \subseteq K$}
		\State Create an arc $a_{JK}$ from $J$ to $K$
		\State $I(a_{JK}) \gets K \setminus J$
		\State $\hat\setA \gets \hat\setA \cup \{a_{JK}\}$
		\EndIf
		\EndFor
		\EndFor

		\State \Return ($\hat\setV$, $\hat\setA$)
	\end{algorithmic}
\end{algorithm}

The single-level reformulation (\cref{prog:dp-single}) derived from the initial selection diagram
is a relaxation of the one derived from the full selection diagram.
Indeed, we did not connect any node to the terminal node $q$.
As a result, the value of $y_\varnothing$ is not ``grounded'', \ie it is in no constraint involving other $y$ variables, and the dual representation has no effect at the beginning.
Over the course of the procedure, we introduce new paths by connecting nonterminal nodes $J$ to $q$
using a \textit{modified terminal arc} $\ha_J$.
As long as the paths that we add to the diagram correspond to feasible solutions in $\setX$,
the reformulation remains valid. Such solutions are obtained by solving the follower's problem similarly to \cref{alg:vf-dynamic}.
The whole process is described in \cref{alg:sdiag-dynamic}, where the difference from \cref{alg:vf-dynamic} is that lines 6-7 are replaced by lines 6-18. Essentially, this part of the algorithm provides a way to represent any feasible solution on $\mathcal{X}$ as a path in the selection diagram.

\begin{algorithm}[ht]
	\caption{Solution of \cref{prog:dp-single} using selection diagram with dynamically generated constraints}
	\label{alg:sdiag-dynamic}

	\begin{algorithmic}[1]
		\Require{The CPP instance, the initial selection diagram $(\hat\setV, \hat\setA)$.}
		\Ensure{An optimal solution $(t^*, x^*)$ of the CPP instance.}
		\State Create the master problem from \cref{prog:dp-single} using the diagram $(\hat\setV, \hat\setA)$.
		\Loop
		\State Solve the master problem and let $(\tilde{t}, \tilde{x})$ be the current optimal solution.
		\State Solve \cref{prog:cpp-follower} with $t = \tilde{t}$ for $\hx$.
		\Comment{\ie choose any $\hx \in \funcR(\tilde{t})$}

		\If{$\tr{(v-\tilde{t})}\tilde{x} < \tr{(v-\tilde{t})}\hat{x}$}
		\LComment{$(\tilde{t},\tilde{x})$ is not bilevel feasible. Add new constraint.}
		\State $added \gets \mathbf{false}$
		\For{$k = 2,1,0$}
		\ForAll{$J \in \hat\setV_k$} 
		\Comment{Iterate in a random order}
		\If{$J \subseteq \setI(\hx)$}
		\State Create an arc $\ha_J$ from $J$ to $q$
		\State $I(\ha_J) \gets \setI(\hx) \setminus J$
		\State $\hat\setA \gets \hat\setA \cup \{\ha_J\}$
		\State Add $y_J \geq \sum_{i \in I(\ha_J)} (v_i - t_i)$ to the master program.
		\State $added \gets \mathbf{true}$
		\State \textbf{break}
		\EndIf
		\EndFor
		\If{$added$}
		\State \textbf{break}
		\Comment{Break the loop, go back to the outer loop}
		\EndIf
		\EndFor

		\Else
		\LComment{$(\tilde{t},\tilde{x})$ is bilevel feasible. Stop.}
		\State $(t^*, x^*) \gets (\tilde{t}, \tilde{x})$
		\State \textbf{break}
		\EndIf
		\EndLoop
		\State \Return $(t^*, x^*)$
	\end{algorithmic}
\end{algorithm}

\begin{example}
	\label{ex:kpp-sd-dynamic}
	Consider the KPP instance introduced in \cref{ex:kpp-vf}.
	At line 4 in \cref{alg:pair-sampling},
	instead of sampling a random solution, we once again exploit monotonicity (\cref{rm:kpp-maximal})
	and sample only maximal solutions in $\setX$.
	Suppose that $N=2$ and that we sample the two pairs $\{1,2\}$ and $\{2,4\}$ in line 5.
	The initial selection diagram returned by \cref{alg:pair-sampling} is drawn in \cref{fig:kpp-sd-dynamic}
	(consider only the solid arcs).

	The dashed arcs in \cref{fig:kpp-sd-dynamic} represent the modified terminal arcs added by \cref{alg:sdiag-dynamic}
	corresponding to the four maximal solutions.
	The two solutions $\{1,2,3\}$ and $\{2,4\}$ contain some pairs in the initial diagram, hence we can connect those pairs
	to $q$ (which is the best case).
	There is no pair contained in the other two solutions $\{1,4\}$ and $\{3,4\}$,
	so we fall back to singletons.
	Note that some solutions like $\{1,4\}$ can be introduced in multiple ways (from either state $\{1\}$ or $\{4\}$).
	In that case, we choose a random node by iterating $\hat\setV_k$ in a random order (line 9 in \cref{alg:sdiag-dynamic}).
	In the worst case, line~10 in \cref{alg:sdiag-dynamic} will be true for $k=0$ since $\hat\setV_0 = \{\varnothing\}$, and we can add an arc from $\varnothing$ directly to $q$, which is equivalent to
	a value function constraint. Thus, we can always introduce any feasible solution into the diagram.

	\begin{figure}[ht]
		\centering
		\scalebox{0.9}{\begin{tikzpicture}[xscale=2, yscale=1.2]
	\tikzstyle{state}=[draw,circle,inner sep=0pt,minimum size=24pt]
	\tikzstyle{action}=[circle,fill=white,inner sep=1pt]

	\node[state] (S0) at (0, -1) { $\varnothing$ };

	\node[state] (S1) at (1, -0) { $\{1\}$ };
	\node[state] (S2) at (1, -1) { $\{2\}$ };
	\node[state] (S4) at (1, -2) { $\{4\}$ };

	\node[state] (S12) at (2, -0.5) { \small $\{1,2\}$ };
	\node[state] (S24) at (2, -1.5) { \small $\{2,4\}$ };

	\node[state] (BS) at (3, -1) { $q$ };

    \draw (S0) -- (S1);
    \draw (S0) -- (S2);
    \draw (S0) -- (S4);

    \draw (S1) -- (S12);
    \draw (S2) -- (S12);
    \draw (S2) -- (S24);
    \draw (S4) -- (S24);

    \draw[dashed] (S12) -- node[action] {$\{3\}$} (BS);
    \draw[dashed] (S24) -- node[action] {$\varnothing$} (BS);
    \draw[dashed] (S1) to[bend left=60] node[action] {$\{4\}$} (BS);
    \draw[dashed] (S4) to[bend right=60] node[action] {$\{3\}$} (BS);

\end{tikzpicture}}
		\caption{Dynamically generated selection diagram of \cref{ex:kpp-sd-dynamic}.
			Solid arcs are the initial arcs, while dashed arcs represent the modified terminal arcs added during \cref{alg:sdiag-dynamic}.
			The labels of the dashed arcs are $I(\ha_J)$.}
		\label{fig:kpp-sd-dynamic}
	\end{figure}
\end{example}

\subsection{Decision Diagram}
\label{ssec:ddiag}

A family of dynamic programming models that is commonly used in the literature is \textit{decision diagrams} \citep{bergman2016,lozano2022}.
In a decision diagram, each step represents the decision to include or exclude a specific item.
Let $n = |\setI|$ and let $(i_{1}, i_{2}, \ldots, i_{n})$ be an ordering of the items in $\setI$.
A \textit{proper node} (neither $p$ nor $q$) in a decision diagram is a tuple $(k, \sk)$ of a step $k \in \{1,\ldots,n-1\}$ and a state $\sk$.
The definition of the states $\sk$ is problem-specific.
We partition $\setV$ into $n + 1$ layers:
$\setV_0$ which contains only $p$; $\setV_k$ which contains all nodes $(k, \sk)$ for $k=1,\ldots,n-1$;
and $\setV_n$ which contains only $q$.

For $k = 0,\ldots,n-1$, a node $u_k \in \setV_k$ has at most 2 outgoing arcs $a^0_{u_k}$ and $a^1_{u_k}$,
corresponding to the decisions $x_{i_{k+1}} = 0$ and $x_{i_{k+1}} = 1$.
The targets of these arcs are in the next layer $\setV_{k+1}$.
The item function is defined as follows:
\begin{align*}
	I(a^0_{u_k}) & = \varnothing, &  I(a^1_{u_k}) & = \{i_{k+1}\}.
\end{align*}

A path in a decision diagram represents a series of decisions $(x_{i_1}, \ldots, x_{i_n})$,
fixing the choice of one item at every step.
At the end of the path, all entries in $x$ have been fixed and the value of $x$ should
correspond to a feasible solution $\hx \in \setX$ given a suitable description of the nodes.
In this paper, we provide the descriptions of the decision diagrams for 3 specific problems:
the knapsack problem (\cref{ex:kpp-dd}), the maximum stable set problem (\cref{sssec:dd-maxstab}), and the minimum set cover problem (\cref{sssec:dd-mincover}).
The descriptions of the decision diagrams for several other problems can be found in \citep{bergman2016}.

\begin{example}
	\label{ex:kpp-dd}
	In the decision diagram of the knapsack problem,
	we define the state $\sk$ at step $k$ as the remaining capacity after $x_{i_1},\ldots,x_{i_k}$ have been decided.
	Furthermore, we express the initial node $p$ as $(0, C)$ and the terminal node $q$ as $(n, 0)$.
	The arc $a^0_{u_k}$ is always available and it sends $(k,\sk)$ to $(k+1,\sk)$.
	The arc $a^1_{u_k}$ only exists if $\sk \geq w_{i_{k+1}}$,
	\ie there is enough remaining capacity to fit the item $i_{k+1}$.
	In that case, it connects $(k,\sk)$ to $(k+1,\sk - w_{i_{k+1}})$.
	An exception is that all arcs from the second-to-last layer $\setV_{n-1}$ end at the terminal node instead.

	The decision diagram for the KPP instance in \cref{ex:kpp-vf} is illustrated in \cref{fig:kpp-dd-full}.
	We can use \cref{rm:kpp-maximal} to simplify the graph.
	At step $k$, if the state $\sk$ satisfies $\sk \geq \sum_{j=k+1}^{n} w_{i_j}$,
	then we must set $x_{i_j} = 1$ for all $j > k$ in order to have a maximal solution.
	Thus, for each $k$, we can merge all states $\sk$ such that $\sk \geq \sum_{j=k+1}^{n} w_{i_j}$.
	Furthermore, if $a^0_{u_k}$ and $a^1_{u_k}$ end at the same node, we can eliminate $a^0_{u_k}$ since including the item $i_k$ is prioritized\footnotemark. 	The final result is shown in \cref{fig:kpp-dd-simplified}.
	Applying \cref{prog:dp-single} to \cref{fig:kpp-dd-simplified},
	\cref{cons:dp-single-arc,cons:dp-single-term} become (note that $y_p = y_{0,3}$):
	\begin{align*}
		y_{0,3} & \geq y_{1,3}, & y_{0,3} & \geq 1 - t_{1} + y_{1,2}, & y_{2,1} & \geq y_{3,1}, & y_{2,1} & \geq 1 - t_{3} + y_{3,0}, \\
		y_{1,2} & \geq y_{2,2}, & y_{1,2} & \geq 1 - t_{2} + y_{2,1}, & y_{2,2} & \geq y_{3,2}, & y_{2,2} & \geq 1 - t_{3} + y_{3,1}, \\
		y_{1,3} & \geq y_{2,3}, & y_{1,3} & \geq 1 - t_{2} + y_{2,2}, & y_{3,0} & \geq 0,       & y_{2,3} & \geq 1 - t_{3} + y_{3,2}, \\
		        &               &         &                           & y_{3,1} & \geq 0,       & y_{3,2} & \geq 1.
	\end{align*}
  \footnotetext{This can happen, for example, when $w_{k}=0$ or at the destination node after merge.}

	\begin{figure}[ht]
		\centering
		\subfloat[Full]{
			\label{fig:kpp-dd-full}
			\scalebox{0.9}{\begin{tikzpicture}[xscale=1.3]
	\tikzstyle{state}=[draw,circle,inner sep=0pt,minimum size=20pt]
	\tikzstyle{act1}=[solid]
	\tikzstyle{act0}=[dashed]

	\node[state] (S13) at (0, 0) { 3 };

	\node[state] (S22) at (1, 1) { 2 };
	\node[state] (S23) at (1, 0) { 3 };

	\node[state] (S31) at (2, 2) { 1 };
	\node[state] (S32) at (2, 1) { 2 };
	\node[state] (S33) at (2, 0) { 3 };

	\node[state] (S40) at (3, 3) { 0 };
	\node[state] (S41) at (3, 2) { 1 };
	\node[state] (S42) at (3, 1) { 2 };
	\node[state] (S43) at (3, 0) { 3 };

	\node[state] (S50) at (4, 3) { 0 };

  \draw[act0] (S40) -- (S50);
  \draw[act0] (S31) -- (S41) -- (S50);
  \draw[act0] (S22) -- (S32) -- (S42);
  \draw[act0] (S13) -- (S23) -- (S33) -- (S43);

  \draw[act1] (S13) -- (S22) -- (S31) -- (S40);
  \draw[act1] (S23) -- (S32) -- (S41);
  \draw[act1] (S33) -- (S42);
  \draw[act1] (S43);

  \draw[decoration={dashsoliddouble}, decorate] (S50) -- (S42);
  \draw[decoration={dashsoliddouble}, decorate] (S50) -- (S43);

\end{tikzpicture}}
		}
		\quad
		\subfloat[Simplified]{
			\label{fig:kpp-dd-simplified}
			\scalebox{0.9}{\begin{tikzpicture}[xscale=1.3]
	\tikzstyle{state}=[draw,circle,inner sep=0pt,minimum size=20pt]
	\tikzstyle{act1}=[solid]
	\tikzstyle{act0}=[dashed]

	\node[state] (S13) at (0, 0) { 3 };

	\node[state] (S22) at (1, 1) { 2 };
	\node[state] (S23) at (1, 0) { 3 };

	\node[state] (S31) at (2, 2) { 1 };
	\node[state] (S32) at (2, 1) { 2 };
	\node[state] (S33) at (2, 0) { 3 };

	\node[state] (S40) at (3, 3) { 0 };
	\node[state] (S41) at (3, 2) { 1 };
	\node[state] (S42) at (3, 1) { 2 };

	\node[state] (S50) at (4, 3) { 0 };

  \draw[act0] (S40) -- (S50);
  \draw[act0] (S31) -- (S41) -- (S50);
  \draw[act0] (S22) -- (S32) -- (S42);
  \draw[act0] (S13) -- (S23) -- (S33);

  \draw[act1] (S13) -- (S22) -- (S31) -- (S40);
  \draw[act1] (S23) -- (S32) -- (S41);
  \draw[act1] (S33) -- (S42) -- (S50);
  \draw[act1] (S43);

\end{tikzpicture}}
		}
		\caption{Decision diagrams of \cref{ex:kpp-dd}.
			The labels of the nodes are the states $\sk$. The steps of the nodes are implied from left to right as $k = 0,\ldots,4$.
			Solid arcs and dashed arcs represent the decisions $x_k = 1$ and $x_k = 0$, respectively.}
		\label{fig:kpp-dd}
	\end{figure}
\end{example}

\subsubsection{Dynamic Generation}
\label{ssec:dynamic-gen}
Similar to selection diagrams, the number of nodes in a decision diagram is exponential on $|\setI|$.
Thus, we also wish to solve \cref{prog:dp-single} using a dynamically generated decision diagram.
\Cref{alg:dd-initial} describes a sampling algorithm that returns an initial decision diagram.
The size of the initial decision diagram is controlled by a parameter $W$ which is called the \textit{width} of the diagram.
As the name suggests, it holds that $|\hat\setV_k| \leq W$ for all $k=1,\ldots,n-1$.
The algorithm samples exactly $W$ solutions $\hx \in \setX$ (line 4), finds a path in the full diagram for each $\hx$ (line 5),
then incorporates those paths into the initial diagram (lines 6 and 7).
The path-finding step is problem-specific, but usually, it can be done quite efficiently without explicitly enumerating the full diagram.

\begin{algorithm}[ht]
	\caption{Generation of the initial decision diagram}
	\label{alg:dd-initial}

	\begin{algorithmic}[1]
		\Require{The CPP instance, the desired width \(W\).}
		\Ensure{The set of nodes $\hat\setV$ and the set of initial arcs $\hat\setA$.}
		\State $\hat\setV \gets \varnothing$
		\Comment{Initial nodes}
		\State $\hat\setA \gets \varnothing$
		\Comment{Initial arcs}

		\For{$\omega = 1,\ldots,W$}
		\State Sample $\hx \in \setX$
		\State Find a path $\pi$ corresponding to $\hx$ in the full decision diagram
		\State $\hat\setV \gets \hat\setV \cup \setV(\pi)$
		\Comment{$\setV(\pi)$ is the set of nodes in $\pi$}
		\State $\hat\setA \gets \hat\setA \cup \setA(\pi)$
		\Comment{$\setA(\pi)$ is the set of arcs in $\pi$}
		\EndFor

		\State \Return ($\hat\setV$, $\hat\setA$)
	\end{algorithmic}
\end{algorithm}

Updating a decision diagram is more complex than updating a selection diagram.
We propose a method to append a path corresponding to a new solution $\hx \in \setX$ using as many nodes in $\hat\setV$ as possible.
Let $\setB$ be a superset of $\setA$ consisting of all \textit{valid transitions} of the decision diagram, which is problem-specific.
The longest path in $(\setV, \setB)$ must have the same length as the longest path in $(\setV, \setA)$.
Unlike $\setA$, for $b \in \setB$, $u(b)$ and $w(b)$ need not be in consecutive layers,
and $I(b)$ can have more than one item.
For instance, $\setB$ can include, but is not limited to, the aggregations of path segments in $\setA$.
The testing of membership in $\setB$ should be fast and should not require the enumeration of $\setB$.

We partition $\hat{\setV}$ into layers $\hat{\setV}_0,\ldots,\hat{\setV}_n$ such that $\hat{\setV}_k \subseteq \setV_k$ for all $k = 0,\ldots,n$.
Given a solution $\hx \in \setX$ that we wish to introduce into the diagram, consider the subset:
\begin{equation*}
	\hat\setB(\hx) = \{ b \in \setB \mid u(b) \in \hat{\setV}_j,~w(b) \in \hat{\setV}_k,~I(b) = \setI_{jk}(\hx) ~\text{for some}~ j < k\}
\end{equation*}
where
\begin{equation*}
	\setI_{jk}(\hx) = \setI(\hx) \cap \{i_{j+1},\ldots,i_k\}.
\end{equation*}

Enumerating $\hat{\setB}(\hx)$ takes $\bigO(|\hat\setV|^2)$-time
assuming that the testing of membership in $\setB$ takes constant time
(a loop for $u(b) \in \hat{\setV}_j$ and an inner loop for $w(b) \in \hat{\setV}_k$).
By construction, every path in the graph $(\hat{\setV}, \hat{\setB}(\hx))$ corresponds to $\hx$.
We then find the longest path in $(\hat{\setV}, \hat{\setB}(\hx))$ where the length of all arcs is exactly 1.
Finally, we introduce this path as a set of new constraints.
The whole process is outlined in \cref{alg:ddiag-dynamic}.

\begin{algorithm}[ht]
	\caption{Solution of \cref{prog:dp-single} using decision diagram with dynamically generated constraints}
	\label{alg:ddiag-dynamic}

	\begin{algorithmic}[1]
		\Require{The CPP instance, the initial decision diagram $(\hat\setV, \hat\setA)$.}
		\Ensure{An optimal solution $(t^*, x^*)$ of the CPP instance.}
		\State Create the master problem from \cref{prog:dp-single} using the diagram $(\hat\setV, \hat\setA)$.
		\Loop
		\State Solve the master problem and let $(\tilde{t}, \tilde{x})$ be the current optimal solution.
		\State Solve \cref{prog:cpp-follower} with $t = \tilde{t}$ for $\hx$.
		\Comment{\ie choose any $\hx \in \funcR(\tilde{t})$}

		\If{$\tr{(v-\tilde{t})}\tilde{x} < \tr{(v-\tilde{t})}\hat{x}$}
		\LComment{$(\tilde{t},\tilde{x})$ is not bilevel feasible. Add new constraint.}
		\State Enumerate $\hat{\setB}(\hx)$.
		\State Find the longest path $\pi$ in $(\hat{\setV}, \hat{\setB}(\hx))$ where all arcs have length 1.
		\State $\hat\setA \gets \hat\setA \cup \setA(\pi)$
		\Comment{$\setA(\pi)$ is the set of arcs in $\pi$}
		\For{$a \in \setA(\pi)$}
		\State Add $y_{u(a)} \geq \sum_{i \in I(a)} (v_i - t_i) + y_{w(a)}$ to the master program.
		\EndFor

		\Else
		\LComment{$(\tilde{t},\tilde{x})$ is bilevel feasible. Stop.}
		\State $(t^*, x^*) \gets (\tilde{t}, \tilde{x})$
		\State \textbf{break}
		\EndIf
		\EndLoop
		\State \Return $(t^*, x^*)$
	\end{algorithmic}
\end{algorithm}

\begin{example}\label{ex:kpp-dd-dynamic}
	In the KPP, we can derive a path corresponding to a given $\hx \in \setX$
	in the full diagram (line 5 in \cref{alg:dd-initial}) by keeping track of the remaining capacity at step $k$,
	defined by the following recursive equations:
	\begin{align*}
		s_0 & = C,                                                   \\
		s_k & = \begin{cases}
			        s_{k-1} - w_{i_k} & \text{if}~i_k \in \setI(\hx),    \\
			        s_{k-1}           & \text{if}~i_k \notin \setI(\hx),
		        \end{cases} &  & \text{for}~k = 1,\ldots,n-1.
	\end{align*}
	The nodes on the path are $p,(1,s_1),\ldots,(n-1,s_{n-1}),q$.
	Applying \cref{alg:dd-initial} to the KPP instance in \cref{ex:kpp-vf}
	with two sampled solutions $\{1,2,3\}$ and $\{3,4\}$,
	we have the initial decision diagram drawn in \cref{fig:kpp-dd-dynamic} (consider only solid and dashed arcs).

	We define the set of valid transitions $\setB$ for the KPP to be:
	\begin{equation*}
		\setB = \left\{ b ~\middle|~
		\begin{array}{lll}
			u(b) = (j, s_j) \in \setV_j, & \quad & j = 0,\ldots,n-1; \\
			w(b) = (k, s_k) \in \setV_k, &       & k = j+1,\ldots,n; \\
			s_j - \sum_{i \in I(b)} w_i \geq s_k
		\end{array}\right\}.
	\end{equation*}
	In other words, a transition $b$ from $(j,s_j)$ to $(k,s_k)$ is valid if and only if
	all items in $I(b)$ can be fit in a knapsack with capacity $s_j - s_k$.
	Recall that $p$ and $q$ can be expressed as $(0, C)$ and $(n, 0)$.
	It is evident that given $u(b)$, $w(b)$, and $I(b)$, this condition can be tested efficiently.
	If we follow \cref{alg:ddiag-dynamic} and add the paths corresponding to $\{1,4\}$ and $\{2,4\}$,
	we will obtain the dotted arcs in \cref{fig:kpp-dd-dynamic}.

	\begin{figure}[ht]
		\centering
		\scalebox{0.9}{\begin{tikzpicture}[xscale=2.0, yscale=1.8]
	\tikzstyle{state}=[draw,circle,inner sep=0pt,minimum size=20pt]
    \tikzstyle{action}=[circle,fill=white,inner sep=1pt,font={\small}]
	\tikzstyle{act1}=[solid]
	\tikzstyle{act0}=[dashed]
	\tikzstyle{trans}=[dotted,thick]

	\node[state] (S13) at (0, 0.5) { 3 };

	\node[state] (S22) at (1, 1) { 2 };
	\node[state] (S23) at (1, 0) { 3 };

	\node[state] (S31) at (2, 1) { 1 };
	\node[state] (S33) at (2, 0) { 3 };

	\node[state] (S40) at (3, 1) { 0 };
	\node[state] (S42) at (3, 0) { 2 };

	\node[state] (S50) at (4, 0.5) { 0 };

    \draw[act1] (S13) -- node[action] {$\{1\}$} (S22) -- node[action] {$\{2\}$} (S31) -- node[action] {$\{3\}$} (S40);
    \draw[act0] (S40) -- node[action] {$\varnothing$} (S50);

    \draw[act0] (S13) -- node[action] {$\varnothing$} (S23) -- node[action] {$\varnothing$} (S33);
    \draw[act1] (S33) -- node[action] {$\{3\}$} (S42) -- node[action] {$\{4\}$} (S50);

	\draw[trans] (S22) -- node[action] {$\varnothing$} (S42);
    \draw[trans] (S23) to[bend right=45] node[action] {$\{2\}$} (S42);

\end{tikzpicture}}
		\caption{Dynamically generated decision diagram of \cref{ex:kpp-dd-dynamic}.
			Solid and dashed arcs are the initial arcs obtained from \cref{fig:kpp-dd-simplified},
			while dotted arcs are valid transitions added during \cref{alg:ddiag-dynamic}.
			The labels of the nodes are the states $\sk$.
			The labels of the arcs are $I(a)$.}
		\label{fig:kpp-dd-dynamic}
	\end{figure}
\end{example}

\subsubsection{Item Grouping}
\label{ssec:grouping}

A major disadvantage of decision diagrams compared to selection diagrams is the lack of scalability.
For every new item, we need to add a whole new layer to the decision diagram, while the number of layers in a selection diagram stays at four.
Furthermore, many nodes just repeat the states of some nodes in the previous layers (see \cref{fig:kpp-dd-full}).
In this section, we describe a technique to reduce the number of nodes in a decision diagram
by grouping multiple items into one layer, instead of keeping the item-layer ratio at 1:1.

Let $(J_1,\ldots,J_m)$ be a partition of $\setI$.
We group all items in each \textit{part} $J_k$ together into one layer, so overall, the grouped decision diagram has $m + 1$ layers $\setV_0,\setV_1,\ldots,\setV_m$.
For $k = 0,\ldots,m-1$, each node $u_k \in \setV_k$ can have up to $2^{|J_{k+1}|}$ outgoing arcs towards the next layer $\setV_{k+1}$.
Each arc $a^K_{u_k}$ corresponds to a subset $K$ of $J_{k+1}$ with the item function simply defined as $I(a^K_{u_k}) = K$.
Such an arc represents the decision of fixing $x_i = 1$ for all $i \in K$, and fixing $x_i = 0$ for all $i \in J_{k+1} \setminus K$.
A path in the grouped diagram still represents a solution $\hx \in \setX$
because $\{J_1,\ldots,J_m\}$ forms a partition of $\setI$. Hence, all entries in $x$ are fixed and no entry is fixed twice.
\Cref{alg:dd-initial,alg:ddiag-dynamic} are extended accordingly.

\begin{example}\label{ex:kpp-dd-group}
	Consider the decision diagram in \cref{fig:kpp-dd-simplified}.
	We group the first two layers into one group $J_1 = \{1,2\}$ and the last two layers into another group $J_2 = \{3,4\}$.
	The result is shown in \cref{fig:kpp-dd-grouped-full}. The grouped version of \cref{fig:kpp-dd-dynamic}
	is illustrated in \cref{fig:kpp-dd-grouped-dyn}.

	\begin{figure}[ht]
		\centering
		\subfloat[Full diagram]{
			\label{fig:kpp-dd-grouped-full}
			\scalebox{0.9}{\begin{tikzpicture}[xscale=1.3, yscale=1.2]
	\tikzstyle{state}=[draw,circle,inner sep=0pt,minimum size=20pt]
    \tikzstyle{action}=[rectangle,fill=white,inner sep=1pt,font={\small}]
	\tikzstyle{act1}=[solid]
	\tikzstyle{act0}=[dashed]
	\tikzstyle{trans}=[dotted,thick]

	\node[state] (S13) at (0, 1) { 3 };

	\node[state] (S31) at (2, 2) { 1 };
	\node[state] (S32) at (2, 1) { 2 };
	\node[state] (S33) at (2, 0) { 3 };

	\node[state] (S50) at (4, 1) { 0 };

    \draw[act0] (S13) to[bend right=20] node[action] {$\varnothing$} (S33);
    \draw[act1] (S13) to[bend right=20] node[action] {$\{1\}$} (S32);
    \draw[act1] (S13) to[bend left=20] node[action] {$\{2\}$} (S32);
    \draw[act1] (S13) to[bend left=20] node[action] {$\{1,2\}$} (S31);

    \draw[act1] (S33) to[bend right=20] node[action] {$\{3,4\}$} (S50);
    \draw[act1] (S32) to[bend right=20] node[action] {$\{4\}$} (S50);
    \draw[act1] (S32) to[bend left=20] node[action] {$\{3\}$} (S50);
    \draw[act1] (S31) to[bend left=20] node[action] {$\{3\}$} (S50);

\end{tikzpicture}}
		}
		\qquad
		\subfloat[Dynamically generated]{
			\label{fig:kpp-dd-grouped-dyn}
			\scalebox{0.9}{\begin{tikzpicture}[xscale=1.3, yscale=1.2]
	\tikzstyle{state}=[draw,circle,inner sep=0pt,minimum size=20pt]
    \tikzstyle{action}=[rectangle,fill=white,inner sep=1pt,font={\small}]
	\tikzstyle{act1}=[solid]
	\tikzstyle{act0}=[dashed]
	\tikzstyle{trans}=[dotted,thick]

	\node[state] (S13) at (0, 1) { 3 };

	\node[state] (S31) at (2, 2) { 1 };
	\node[state] (S33) at (2, 0) { 3 };

	\node[state] (S50) at (4, 1) { 0 };

    \draw[act0] (S13) -- node[action] {$\varnothing$} (S33);
    \draw[act1] (S13) -- node[action] {$\{1,2\}$} (S31);

    \draw[act1] (S33) -- node[action] {$\{3,4\}$} (S50);
    \draw[act1] (S31) -- node[action] {$\{3\}$} (S50);

    \draw[trans] (S13) to[bend right=10] node[action] {$\{1,4\}$} (S50);
    \draw[trans] (S13) to[bend left=10] node[action] {$\{2,4\}$} (S50);

\end{tikzpicture}}
		}
		\caption{Grouped decision diagrams of \cref{ex:kpp-dd-group}.
			Solid and dashed arcs are arcs in the initial diagram,
			dotted arcs are generated valid transitions.
			The labels of the nodes are the states $\sk$.
			The labels of the arcs are $I(a)$.}
		\label{fig:kpp-dd-grouped}
	\end{figure}
\end{example}

In our experiments, we keep the number of layers $m$ fixed and we partition the items randomly into groups of size $\lceil n/m \rceil$.
This strategy allows us to maintain a constant number of nodes $|\hat{\setV}|$ in the initial decision diagram as $n$ increases.
We remark that this technique should be interpreted as a general framework,
and one can develop new strategies depending on the situation, including groups with different sizes,
subgroups inside a group, separate groups for tolled items and toll-free items, etc.

\subsubsection{Maximum Stable Set Pricing Problem}
\label{sssec:dd-maxstab}
We consider the case when the follower's problem is a maximum stable set problem.
The overall pricing problem is called the \textit{maximum stable set pricing problem} (MaxSSPP).
Let $\setG = (\setI, \setE)$ be an undirected simple graph.
In this case, the set of items $\setI$ is the set of nodes of $\setG$.
The objective of the follower is to find a \textit{stable set} (\ie a subset of $\setI$ in which no pair of nodes are adjacent)
that maximizes the sum of tolled values.
If $A \in \setbin^{\setI \times \setE}$ is the node-edge incidence matrix of $\setG$,
then the follower's problem of the MaxSSPP is formulated as follows:
\begin{equation*}
	\funcR(t) = \underset{x}{\arg\max} \{ \vtx \mid \tr{A}x \leq \vecone,~x \in \setbin^{\setI} \}
\end{equation*}
where $\vecone$ is the vector of all ones.

In the decision diagram for the MaxSSPP, the state $\sk$ at step $k$
is the set of all items in $\{i_{k+1},\ldots,i_n\}$ that are still available,
\ie they are not adjacent to any already-selected item.
The arc $a^0_{u_k}$ is always available and it sends $(k, \sk)$ to $(k+1, \sk \setminus \{i_{k+1}\})$.
The arc $a^1_{u_k}$ only exists if $i_{k+1} \in \sk$ (meaning the item $i_{k+1}$ is available).
In that case, it sends $(k, \sk)$ to $(k+1, \sk \setminus N[i_{k+1}])$
where $N[i]$ is the closed neighborhood of $i$ (\ie the set that includes $i$ and all nodes adjacent to $i$).
We express $p$ as $(0, \setI)$ and $q$ as $(n, \varnothing)$.

In regard to \cref{alg:dd-initial}, given $\hx \in \setX$,
we can generate a path corresponding to $\hx$ by the following recursive system:
\begin{align*}
	s_0 & = \setI,                                                       \\
	s_k & = \begin{cases}
		        s_{k-1} \setminus N[i_k]  & \text{if}~i_k \in \setI(\hx),    \\
		        s_{k-1} \setminus \{i_k\} & \text{if}~i_k \notin \setI(\hx),
	        \end{cases} &  & \text{for}~k = 1,\ldots,n.
\end{align*}

The last piece is the definition of $\setB$ used in \cref{alg:ddiag-dynamic}.
A transition $b$ from $(j, s_j)$ to $(k, s_k)$ is valid if and only if
$I(b)$ is a stable set, $I(b) \subseteq s_j$, and $s_k \subseteq s_j \setminus N[I(b)]$.

\begin{example}\label{ex:dd-maxstab}
	Consider a MaxSSPP instance of 5 items whose graph is shown in \cref{fig:maxstab-instance}.
	Assuming that the item ordering is $(1,2,3,4,5)$,
	the initial decision diagram obtained from \cref{alg:dd-initial} with two sampled solutions
	$\{1,3\}$ and $\{2,5\}$ is drawn in \cref{fig:maxstab-dd-dynamic} (consider only solid and dashed arcs).
	The dotted arcs in \cref{fig:maxstab-dd-dynamic} represent the valid transitions generated by \cref{alg:ddiag-dynamic}
	corresponding to the two solutions $\{1,4\}$ and $\{3,5\}$.

	\begin{figure}[ht]
		\centering
		\scalebox{0.9}{\begin{tikzpicture}[scale=1.5]
	\tikzstyle{state}=[draw,circle,inner sep=0pt,minimum size=20pt]

    \node[state] (S1) at (0,1) { $1$ };
    \node[state] (S2) at (1,1) { $2$ };
    \node[state] (S3) at (2,0.5) { $3$ };
    \node[state] (S4) at (1,0) { $4$ };
    \node[state] (S5) at (0,0) { $5$ };

    \draw (S1) -- (S2) -- (S3) -- (S4) -- (S5) -- (S1);
    \draw (S2) -- (S4);
\end{tikzpicture}}
		\caption{The graph $\setG$ of the MaxSSPP instance in \cref{ex:dd-maxstab}.}
		\label{fig:maxstab-instance}
	\end{figure}

	\begin{figure}[ht]
		\centering
		\scalebox{0.9}{\begin{tikzpicture}[xscale=2.0, yscale=2.0]
    \arraycolsep=0pt 
	\tikzstyle{state}=[draw,circle,inner sep=0pt,minimum size=24pt]
    \tikzstyle{action}=[rectangle,fill=white,inner sep=1pt,font={\small}]
	\tikzstyle{act1}=[solid]
	\tikzstyle{act0}=[dashed]
	\tikzstyle{trans}=[dotted,thick]

    \node[state] (S0-12345) at (0, 0.5) { \scriptsize $\left\{\begin{array}{c}1,2,\\3,4,5\end{array}\right\}$ };

	\node[state] (S1-34) at (1, 1) { $\{3,4\}$ };
	\node[state] (S1-2345) at (1, 0) { \scriptsize $\left\{\begin{array}{c}2,3,\\4,5\end{array}\right\}$ };

	\node[state] (S2-34) at (2, 1) { $\{3,4\}$ };
	\node[state] (S2-5) at (2, 0) { $\{5\}$ };

	\node[state] (S3-0) at (3, 1) { $\varnothing$ };
	\node[state] (S3-5) at (3, 0) { $\{5\}$ };

	\node[state] (S4-0) at (4, 1) { $\varnothing$ };
	\node[state] (S4-5) at (4, 0) { $\{5\}$ };

	\node[state] (S5-0) at (5, 0.5) { $\varnothing$ };

    \draw[act1] (S0-12345) -- node[action] {$\{1\}$} (S1-34) (S2-34) -- node[action] {$\{3\}$} (S3-0);
    \draw[act0] (S1-34) -- node[action] {$\varnothing$} (S2-34) (S3-0) -- node[action] {$\varnothing$} (S4-0) -- node[action] {$\varnothing$} (S5-0);

    \draw[act1] (S1-2345) -- node[action] {$\{2\}$} (S2-5) (S4-5) -- node[action] {$\{5\}$} (S5-0);
    \draw[act0] (S0-12345) -- node[action] {$\varnothing$} (S1-2345) (S2-5) -- node[action] {$\varnothing$} (S3-5) -- node[action] {$\varnothing$} (S4-5);

    \draw[trans] (S2-34) to[bend right=35] node[action] {$\{4\}$} (S4-0);
    \draw[trans] (S1-2345) to[bend left=35] node[action] {$\{3\}$} (S3-5);

\end{tikzpicture}}
		\caption{Dynamically generated decision diagram of \cref{ex:dd-maxstab}.
		Solid and dashed arcs are the initial arcs $a^1_{u_k}$ and $a^0_{u_k}$ respectively,
		while dotted arcs are generated valid transitions $b$.
		The labels of the nodes are the states $\sk$.
		The labels of the arcs are $I(a)$.}
		\label{fig:maxstab-dd-dynamic}
	\end{figure}
\end{example}

\subsubsection{Minimum Set Cover Pricing Problem}
\label{sssec:dd-mincover}
The last pricing problem that we investigate is the \textit{minimum set cover pricing problem} (MinSCPP),
in which the follower solves a minimum set cover problem.
Let $\setE$ be the set of elements that we wish to cover and
let $\setI$ be a family of subsets of $\setE$ that covers $\setE$ (\ie the union of all sets in $\setI$ is equal to $\setE$).
For the sake of consistency, we refer to subsets in $\setI$ as \textit{items}.
The follower aims to find a subset of $\setI$ that minimizes the sum of tolled values while still covering $\setE$:
\begin{equation*}
	\funcR(t) = \underset{x}{\arg\min} \{ \tr{(v + t)}x \mid Ax \geq \vecone,~x \in \setbin^{\setI} \}
\end{equation*}
where $A \in \setbin^{\setE \times \setI}$ is the incidence matrix between $\setE$ and $\setI$.
We make the assumption that the set of toll-free items $\setI_2$ must cover $\setE$,
otherwise the leader can impose an infinite toll on all tolled items and the pricing problem is unbounded.

The state $\sk$ in the MinSCPP is defined to be the set of uncovered elements of $\setE$.
In contrast to the first two problems, both arcs $a^0_{u_k}$ and $a^1_{u_k}$ are always available.
The former sends $(k, \sk)$ to $(k + 1, \sk)$, and the latter sends $(k, \sk)$ to $(k+1, \sk \setminus i_{k+1})$
(recall that $i_{k+1}$ is a subset of $\setE$).
An exception exists in the next-to-last layer $\setV_{n-1}$ to check if the selected items form a cover:
if $s_{n-1} \neq \varnothing$ (resp. $s_{n-1} \setminus i_n \neq \varnothing$),
then the arc $u^0_{u_{n-1}}$ (resp. $u^1_{u_{n-1}}$) is not available.
The nodes $p$ and $q$ are expressed as $(0, \setE)$ and $(n, \varnothing)$.

Given $\hx \in \setX$, a path corresponding to $\hx$ is generated using the recursive system:
\begin{align*}
	s_0 & = \setE,                                                   \\
	s_k & = \begin{cases}
		        s_{k-1} \setminus i_k & \text{if}~i_k \in \setI(\hx),    \\
		        s_{k-1}               & \text{if}~i_k \notin \setI(\hx),
	        \end{cases} &  & \text{for}~k = 1,\ldots,n.
\end{align*}

Finally, a transition $b$ from $(j, s_j)$ to $(k, s_k)$ is valid if and only if
\begin{equation*}
	s_j \setminus \Big( \bigcup_{i \in I(b)} i \Big) \subseteq s_k.
\end{equation*}

\begin{example}\label{ex:dd-mincover}
	Consider a MinSCPP instance where $\setE = \{a,b,c,d\}$ and $\setI$ consists of
	five sets $i_1 = \{a,b\}$, $i_2 = \{a,c\}$, $i_3 = \{a,d\}$, $i_4 = \{b,c,d\}$, and $i_5 = \{c\}$.
	We order the items according to $(i_1,i_2,i_3,i_4,i_5)$.
	If we sample the solutions $\{i_1,i_3,i_5\}$ and $\{i_2,i_4\}$ for \cref{alg:dd-initial}
	and then generate valid transitions for $\{i_1,i_2,i_3\}$ and $\{i_3,i_4\}$ in \cref{alg:ddiag-dynamic},
	we will obtain the decision diagram illustrated in \cref{fig:mincover-dd-dynamic}.

	\begin{figure}[ht]
		\centering
		\scalebox{0.9}{\begin{tikzpicture}[xscale=2.0, yscale=2.0]
    \arraycolsep=0pt 
	\tikzstyle{state}=[draw,circle,inner sep=0pt,minimum size=24pt]
    \tikzstyle{action}=[rectangle,fill=white,inner sep=1pt,font={\small}]
	\tikzstyle{act1}=[solid]
	\tikzstyle{act0}=[dashed]
	\tikzstyle{trans}=[dotted,thick]

    \node[state] (S0-abcd) at (0, 0.5) { \scriptsize $\left\{\begin{array}{c}a,b,\\c,d\end{array}\right\}$ };

	\node[state] (S1-cd) at (1, 1) { $\{c,d\}$ };
	\node[state] (S1-abcd) at (1, 0) { \scriptsize $\left\{\begin{array}{c}a,b,\\c,d\end{array}\right\}$ };

	\node[state] (S2-cd) at (2, 1) { $\{c,d\}$ };
	\node[state] (S2-bd) at (2, 0) { $\{b,d\}$ };

	\node[state] (S3-c) at (3, 1) { $\{c\}$ };
	\node[state] (S3-bd) at (3, 0) { $\{b,d\}$ };

	\node[state] (S4-c) at (4, 1) { $\{c\}$ };
	\node[state] (S4-0) at (4, 0) { $\varnothing$ };

	\node[state] (S5-0) at (5, 0.5) { $\varnothing$ };

    \draw[act1] (S0-abcd) -- node[action] {$\{i_1\}$} (S1-cd) (S2-cd) -- node[action] {$\{i_3\}$} (S3-c) (S4-c) -- node[action] {$\{i_5\}$} (S5-0);
    \draw[act0] (S1-cd) -- node[action] {$\varnothing$} (S2-cd) (S3-c) -- node[action] {$\varnothing$} (S4-c);

    \draw[act1] (S1-abcd) -- node[action] {$\{i_2\}$} (S2-bd) (S3-bd) -- node[action] {$\{i_4\}$} (S4-0);
    \draw[act0] (S0-abcd) -- node[action] {$\varnothing$} (S1-abcd) (S2-bd) -- node[action] {$\varnothing$} (S3-bd) (S4-0) -- node[action] {$\varnothing$} (S5-0);

    \draw[trans] (S1-cd) to node[action] {$\{i_2,i_3\}$} (S4-0);
    \draw[trans] (S1-abcd) to[bend right=25] node[action] {$\{i_3,i_4\}$} (S4-0);

\end{tikzpicture}}
		\caption{Dynamically generated decision diagram of \cref{ex:dd-mincover}.
		Solid and dashed arcs are the initial arcs $a^1_{u_k}$ and $a^0_{u_k}$ respectively,
		while dotted arcs are generated valid transitions $b$.
		The labels of the nodes are the states $\sk$.
		The labels of the arcs are $I(a)$.}
		\label{fig:mincover-dd-dynamic}
	\end{figure}
\end{example}

\section{Experiments}
\label{sec:exp}
We conducted numerical experiments to compare the performance of the three dynamic programming models:
the value function reformulation (VF), the selection diagram (SD), and the decision diagram (DD).
All three models were tested on a family of randomly generated instances belonging to three pricing problems:
the knapsack pricing problem (KPP), the maximum stable set pricing problem (MaxSSPP),
and the minimum set cover pricing problem (MinSCPP).
Additionally, we applied these models to a closely related problem that is the
knapsack interdiction problem (KIP).

To evaluate the gradual impact of the selection diagram and the decision diagram compared to the control case (which is the value function reformulation),
we vary the parameters $N$ in \cref{alg:pair-sampling} and $W$ in \cref{alg:dd-initial}.
Note that if $N = 0$ and $W = 0$, then the programs corresponding to both diagrams become the value function reformulation.
To reduce the variance when comparing the results with different values of $N$ and $W$,
we use the same random seed for all sampling procedures (lines 4, 5 in \cref{alg:pair-sampling}, and line 4 in \cref{alg:dd-initial}).
For the decision diagram, we also apply item grouping (\cref{ssec:grouping}) to all instances so that the number of layers is at most 20,
except for the runs on the KPP where we group every 2 items together
(since the number of items in the KPP instances are too low to be evenly divided into 20 layers). We refer the reader to~\cref{app:eff-grouping}, which empirically demonstrates that item grouping is typically beneficial.

Every run was executed on a single thread (AMD EPYC 7532 2.4GHz CPU, 8GB of RAM\footnotemark) with the time limit of one hour.
\footnotetext{The experiment was run on the Narval computing cluster provided by the Digital Research Alliance of Canada.}
Mixed-integer programs are solved using Gurobi 9.5 \citep{gurobi} while the remaining code (formulation and constraint generation) is written in Julia 1.8.5.
The dynamically generated constraints in \cref{alg:vf-dynamic,alg:sdiag-dynamic,alg:ddiag-dynamic}
are implemented as lazy constraints in Gurobi.
We reuse the bilevel feasible solutions $(\tilde{t}, \hx)$ obtained during constraint generation
as feasible solutions of the overall pricing problem.
We supplement the solver with these solutions via a custom heuristic.
Bounds for the McCormick envelope are provided in \cref{app:bounds}.
All instances, results and code are publicly available~\citep{BuiCode2025}. 

The experimental results of each pricing problem are presented in \cref{ssec:exp-kpp,ssec:exp-maxstab,ssec:exp-mincover} respectively. In \cref{ssec:exp-kip}, we describe a modified formulation of \cref{prog:dp-single}
to apply our methodology to the KIP, followed by the corresponding experimental results.
Finally, we provide a summary and our observations in \cref{ssec:discussions}.

\subsection{Knapsack Pricing Problem}
\label{ssec:exp-kpp}
We generated random KPP instances that satisfy the following properties:
\begin{itemize}
	\item The number of items $|\setI|$ ranges from $40$ to $60$.
	\item The weights $w_i$ are uniformly sampled in the range $[1,100]$.
	\item The densities $v_i/w_i$ are sampled uniformly in the range $[0.75, 1.25]$.
	      The values $v_i$ are then calculated by multiplying with the random weights.
	\item The proportion of tolled items $|\setI_1|/|\setI|$ is equal to the ratio between the capacity and the sum of weights, \ie $C/\sum_{i \in \setI}w_i$.
	      We denote this proportion as $r$. Desirable values for $r$ are from $0.5$ to $0.6$.
       The difficulty of the generated instances decreases sharply when $r$ is outside of this range.
	\item The values $v_i$ of the tolled items are multiplied by $2.0$ afterward to encourage the use of tolled items.
        This technique is adopted from \citep{brotcorne2000}.
\end{itemize}
Parameters $C$, $v_i$ and $w_i$ were rounded to the nearest integer. The parameters listed above are tuned using an evolutionary algorithm (CMA-ES \citep{hansen2003})
to maximize the difficulty of the instances (see \cref{app:sens-anal}).
In total, 330 KPP instances were generated according to 33 different configurations (10 instances per configuration).
These 33 configurations are the combinations of eleven  values of $|\setI| \in \{40, 42, \ldots, 58, 60\}$,
and three values of $r \in \{0.50, 0.55, 0.60\}$.

\begin{table}[htbp]
	\centering
    \caption{Numerical results for the KPP.}
    \label{tab:kpp-results}
    \vskip 2mm
	\begin{talltblr}[
        label = none,
		note{a} = {Geometric average.},
		note{b} = {Arithmetic average.},
		note{} = {The best configuration of each model is highlighted in bold.},
		]{
		colspec = {lrrrrr},
		cells = {font=\small},
		cell{1}{-} = {c},
		cell{1}{1,2,3,6} = {r=2}{},
		cell{1}{4} = {c=2}{},
		row{3,6,16} = {font=\bf\small},
			}
		\toprule
		{DP                                      \\ Model} & $N$/$W$ & {Num. of \\ instances \\ solved} & Time\TblrNote{a}~~(s) & & {Gap\TblrNote{b} \\ (\%)} \\
		\cmidrule[lr]{4-5}
		   &      &     & Total & Callback &     \\
		\midrule
		VF &      & 134 & 1301  & 6        & 6.5 \\
		\midrule
		SD & 100  & 118 & 1316  & 10       & 7.6 \\
		SD & 200  & 121 & 1263  & 10       & 7.7 \\
		SD & 300  & 124 & 1263  & 10       & 7.6 \\
		SD & 500  & 120 & 1333  & 9        & 7.7 \\
		SD & 700  & 117 & 1379  & 10       & 7.8 \\
		SD & 1000 & 113 & 1402  & 10       & 8.1 \\
		SD & 1500 & 118 & 1465  & 10       & 8.0 \\
		SD & 2000 & 116 & 1431  & 10       & 7.9 \\
		SD & 2500 & 113 & 1452  & 10       & 8.0 \\
		SD & 3000 & 114 & 1487  & 11       & 8.2 \\
		\midrule
		DD & 10   & 147 & 973   & 5        & 5.7 \\
		DD & 20   & 157 & 994   & 5        & 5.5 \\
		DD & 30   & 158 & 980   & 6        & 5.4 \\
		DD & 50   & 150 & 1082  & 10       & 5.6 \\
		DD & 70   & 144 & 1155  & 14       & 5.6 \\
		DD & 100  & 135 & 1327  & 19       & 6.0 \\
		DD & 150  & 124 & 1540  & 30       & 6.6 \\
		DD & 200  & 113 & 1735  & 42       & 7.1 \\
		DD & 250  & 106 & 1923  & 53       & 7.5 \\
		DD & 300  & 100 & 2070  & 64       & 7.9 \\
		\bottomrule
	\end{talltblr}
\end{table}

\Cref{tab:kpp-results} shows the number of instances solved, the average time, and the average optimality gap
of the three DP models with various values of $N$ and $W$.
The results of the best configurations of each model are plotted in \cref{fig:exp-knapsack-time,fig:exp-knapsack-gap}.
According to \cref{tab:kpp-results}, the selection diagram does not improve the performance of the solution compared to
the value function model, while the decision diagram does improve it slightly but only at low $W$ (under 100). We observe that, compared to the value function model, the selection and decision diagram model do not reduce the callback time. With respect to the subproblems of the next two sections, the knapsack problem is a relatively simple problem, containing only one constraint,
so the time to solve the subproblem (included in the callback time) is negligible. The callback time increases when we increase the width $W$ of the decision diagram,
reflecting the fact that \cref{alg:ddiag-dynamic} scales quadratically with respect to $|\hat{\setV}|$. Therefore, the performance improvement of the decision diagram over the value function is limited by the width $W$. Moreover, the value function model with the dynamically generated constraints is smaller than that of the other two models, which may lead to simpler master programs and could help explain its relatively good performance.

\subsection{Maximum Stable Set Pricing Problem}
\label{ssec:exp-maxstab}
Random MaxSSPP instances were generated on top of random graphs based on the Erd\"os-R\'enyi model.
The parameters for the generation are as follows:
\begin{itemize}
	\item The number of items $|\setI|$ ranges from $120$ to $180$.
	\item The graph density, denoted as $d$, is from $0.12$ to $0.24$.
	\item The proportion of tolled items $|\setI_1|/|\setI|$ is $0.4$.
	\item The values $v_i$ are uniformly sampled in the range $[50, 150]$,
	      then those of the tolled items are multiplied by $1.3$. The obtained values of $v_i$ are then rounded to the nearest integer.
\end{itemize}
Similar to the KPP instances, these parameters are tuned using the same evolutionary process to maximize their difficulty (see~\cref{app:sens-anal}).
We generated 280 instances in total with respect to 28 configurations:
seven values of $|\setI| \in \{120,130,\ldots,180\}$ combined with four values of $d \in \{0.12, 0.16, 0.20, 0.24\}$.

\begin{table}[htbp]
	\centering
    \caption{Numerical results for the MaxSSPP.}
    \label{tab:maxstab-results}
    \vskip 2mm
	\begin{talltblr}[
        label = none,
		note{a} = {Geometric average.},
		note{b} = {Arithmetic average.},
		note{} = {The best configuration of each model is highlighted in bold.},
		]{
		colspec = {lrrrrrrr},
		cells = {font=\small},
		cell{1}{-} = {c},
		cell{1}{1,2,3,6,7,8} = {r=2}{},
		cell{1}{4} = {c=2}{},
		row{3,11,20} = {font=\bf\small},
			}
		\toprule
		{DP                                                    \\ Model} & $N$/$W$ & {Num. of \\ instances \\ solved} & Time\TblrNote{a}~~(s) & &
		{Gap\TblrNote{b}                                       \\ (\%)} & {Num. of \\ callback \\ calls\TblrNote{a}} & {Time per \\ callback \\ call\TblrNote{a}~~(s)} \\
		\cmidrule[lr]{4-5}
		   &      &     & Total & Callback &                   \\
		\midrule
		VF &      & 216 & 654   & 650      & 15.0 & 202 & 3.22 \\
		\midrule
		SD & 100  & 230 & 655   & 650      & 10.0 & 201 & 3.24 \\
		SD & 200  & 228 & 611   & 606      & 10.8 & 194 & 3.12 \\
		SD & 300  & 234 & 610   & 604      & 10.6 & 202 & 3.00 \\
		SD & 500  & 252 & 500   & 493      & 4.9  & 191 & 2.58 \\
		SD & 700  & 250 & 457   & 449      & 3.8  & 185 & 2.42 \\
		SD & 1000 & 257 & 415   & 407      & 3.7  & 181 & 2.25 \\
		SD & 1500 & 254 & 393   & 384      & 4.3  & 174 & 2.20 \\
		SD & 2000 & 265 & 389   & 379      & 2.0  & 172 & 2.21 \\
		SD & 2500 & 261 & 381   & 371      & 2.6  & 168 & 2.20 \\
		SD & 3000 & 259 & 387   & 375      & 2.9  & 167 & 2.25 \\
		\midrule
		DD & 10   & 251 & 543   & 537      & 4.9  & 192 & 2.80 \\
		DD & 20   & 247 & 496   & 489      & 4.4  & 184 & 2.65 \\
		DD & 30   & 250 & 493   & 486      & 4.2  & 188 & 2.58 \\
		DD & 50   & 252 & 502   & 493      & 4.0  & 189 & 2.61 \\
		DD & 70   & 255 & 483   & 474      & 4.2  & 181 & 2.62 \\
		DD & 100  & 253 & 484   & 473      & 4.1  & 186 & 2.54 \\
		DD & 150  & 257 & 487   & 473      & 3.4  & 181 & 2.61 \\
		DD & 200  & 254 & 499   & 481      & 5.1  & 176 & 2.73 \\
		DD & 250  & 258 & 509   & 487      & 3.7  & 175 & 2.78 \\
		DD & 300  & 254 & 478   & 456      & 4.9  & 169 & 2.70 \\
		\bottomrule
	\end{talltblr}
\end{table}

\Cref{tab:maxstab-results,fig:exp-maxstab-time} show the results of the experiments on the MaxSSPP.
The decision diagram still produces a slight improvement and its performance is less sensitive to $W$.
In contrast to the KPP, the selection diagram at high $N$ surpasses the decision diagram.
At $N = 2000$, the selection diagram model can solve 50 more instances and is more than 1.5 times faster compared to the value function model.
Unlike the KPP, the callback time (which includes the time to solve the subproblem) represents a large proportion of
the total solution time. Thus, both diagrams benefit from the decreases in the number of callback calls and the time per callback call.

\subsection{Minimum Set Cover Pricing Problem}
\label{ssec:exp-mincover}
The generation of MinSCPP instances is more complex compared to the other two problems.
First, every element $e \in \setE$ is given its own value $w_e$ sampled uniformly in $[50, 85]$.
Then, we synthesize the sets $i \in \setI$ by including each element in $\setE$ with probability $0.23$.
Special care is required to ensure $\setI$ covers $\setE$.
The value $v_i$ is calculated as the sum of $w_e$ for all $e \in i$, then it is perturbed by a factor sampled within $[0.9, 1.1]$, and rounded to its nearest integer.
Next, we find a subcover of $\setI$ to form the set of toll-free items $\setI_2$
(recall that $\setI_2$ must cover $\setE$ otherwise the leader can impose an infinite toll).
The proportion of tolled items is $0.28$.
Finally, the values of the tolled items are divided by $2.3$.
All parameters are tuned using an evolutionary algorithm, again, to maximize their difficulty (see~\cref{app:sens-anal}).

We experimented with six values of $|\setI| \in \{70,80,\ldots,120\}$ and five values of $|\setI|/|\setE| \in \{0.8, 1.0, 1.2, 1.4, 1.6\}$,
combined to create 30 configurations, hence 300 random MinSCPP instances.

\begin{table}[htbp]
	\centering
    \caption{Numerical results for the MinSCPP.}
    \label{tab:mincover-results}
    \vskip 2mm
	\begin{talltblr}[
        label = none,
		note{a} = {Geometric average.},
		note{b} = {Arithmetic average.},
		note{} = {The best configuration of each model is highlighted in bold.},
		]{
		colspec = {lrrrrrrr},
		cells = {font=\small},
		cell{1}{-} = {c},
		cell{1}{1,2,3,6,7,8} = {r=2}{},
		cell{1}{4} = {c=2}{},
		row{3,10,19} = {font=\bf\small},
			}
		\toprule
		{DP                                                    \\ Model} & $N$/$W$ & {Num. of \\ instances \\ solved} & Time\TblrNote{a}~~(s) & &
		{Gap\TblrNote{b}                                       \\ (\%)} & {Num. of \\ callback \\ calls\TblrNote{a}} & {Time per \\ callback \\ call\TblrNote{a}~~(s)} \\
		\cmidrule[lr]{4-5}
		   &      &     & Total & Callback &                   \\
		\midrule
		VF &      & 267 & 313   & 302      & 10.1 & 305 & 0.99 \\
        \midrule
		SD & 100  & 282 & 259   & 248      & 7.3  & 282 & 0.88 \\
		SD & 200  & 283 & 251   & 240      & 6.5  & 273 & 0.88 \\
		SD & 300  & 283 & 210   & 199      & 5.8  & 238 & 0.84 \\
		SD & 500  & 283 & 130   & 120      & 3.8  & 145 & 0.83 \\
		SD & 700  & 292 & 94    & 84       & 1.0  & 107 & 0.79 \\
		SD & 1000 & 294 & 88    & 79       & 0.4  & 99  & 0.80 \\
		SD & 1500 & 296 & 89    & 80       & 0.4  & 100 & 0.80 \\
		SD & 2000 & 293 & 99    & 89       & 0.3  & 110 & 0.80 \\
		SD & 2500 & 292 & 110   & 99       & 0.5  & 120 & 0.82 \\
		SD & 3000 & 295 & 129   & 116      & 0.3  & 141 & 0.82 \\
        \midrule
		DD & 50   & 282 & 244   & 231      & 7.1  & 264 & 0.88 \\
		DD & 100  & 283 & 232   & 216      & 5.3  & 255 & 0.85 \\
		DD & 150  & 280 & 227   & 208      & 5.7  & 244 & 0.85 \\
		DD & 200  & 284 & 228   & 205      & 5.4  & 241 & 0.85 \\
		DD & 250  & 284 & 221   & 196      & 5.0  & 226 & 0.87 \\
		DD & 300  & 288 & 230   & 203      & 3.5  & 224 & 0.91 \\
		\bottomrule
	\end{talltblr}
\end{table}

The results for the MinSCPP are displayed in \cref{tab:mincover-results,fig:exp-mincover-time}.
Similar to the MaxSSPP, the callback time plays a major role in the total solution time.
The selection diagram greatly outperforms the other two models.
At $N = 1500$, it can solve almost all instances and is 3.5 times faster than the value function model.
The main contribution is the drop in the number of callback calls by three times.
It also reduces the time spent per callback call by 20\%.

\subsection{Knapsack Interdiction Problem}
\label{ssec:exp-kip}

The KIP proposed by~\cite{denegre2011} is not strictly a pricing problem. Nevertheless,
we can utilize the same methodology to derive a single-level
reformulation similar to \cref{prog:dp-single}, showcasing the breadth of our framework, which can be applied to other challenging bilevel programs. For the sake of consistency, let $t,~x \in \{0,1\}^\setI$ denote the leader's decision and the follower's decision, respectively.
Then, the KIP has the following bilevel formulation:
\begin{equation*}
	\min_{t} \{ f(t) \mid \tr{W}t \leq C,~t \in \{0,1\}^\setI \}
\end{equation*}
where $f(t)$ is the optimal objective value of the integer program:
\begin{equation}
    \label[prog]{prog:kip-follower}
	f(t) = \max_{x} \{ \vx \mid \tr{w}x \leq c,~x \leq 1 - t, ~x \in \{0,1\}^\setI \}.
\end{equation}

Due to \citet{caprara2016}, we can rewrite \cref{prog:kip-follower} as:
\begin{equation}
    \label[prog]{prog:kip-bilinear}
	f(t) = \max_{x} \Big\{ \sumI (v_i - v_it_i)x_i \mid \tr{w}x \leq c, ~x \in \{0,1\}^\setI \Big\}.
\end{equation}

We observe that \cref{prog:kip-bilinear} has the same form as \cref{prog:cpp-follower}.
As a result, we can replace \cref{prog:kip-bilinear} with a DP model as in \cref{prog:dp-dual}.
The overall single-level reformulation of the KIP is:
\begin{subequations}
	\label[prog]{prog:kip-single}
	\begin{align}
		\min_{t,x,y}\quad & y_p \label{eq:kip-single-obj}                                                                                 \\
		\st\quad          & \tr{W}t \leq C,                                       \\
                          & \tr{w}x \leq c,                                       \\
		                  & y_{u(a)} \geq F(a; t) + y_{w(a)}, &  &  a \in \setA, \label[cons]{cons:kip-single-arc} \\
		                  & y_q = 0,                          &  & \label[cons]{cons:kip-single-term}             \\
		                  & t,~x \in \{0,1\}^\setI.
	\end{align}
\end{subequations}

Note that the realization of \cref{cons:kip-single-arc,cons:kip-single-term} depends on
the chosen DP model (whether it is the value function, selection diagram, or decision diagram).
Moreover, we reuse $y_p$ in \eqref{eq:kip-single-obj},
since both the leader and the follower in the KIP optimize the same objective function
but in opposite directions. As a side effect, \cref{prog:kip-single} is free from bilinear terms (contrary to the reformulation of \cite{lozano2022})
and hence, it does not require linearization (as described in \cref{ssec:vf-form}).

\begin{figure}[htbp]
	\centering
	\subfloat[KPP - Time (s)]{
		\label{fig:exp-knapsack-time}
		\scalebox{0.85}{\includegraphics{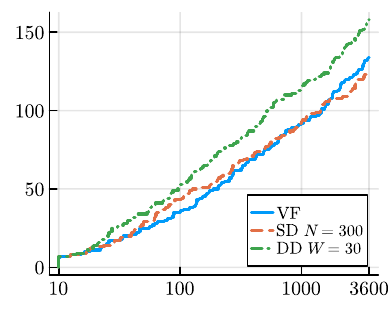}}
	}
	\subfloat[KPP - Optimality gap (\%)]{
		\label{fig:exp-knapsack-gap}
		\scalebox{0.85}{\includegraphics{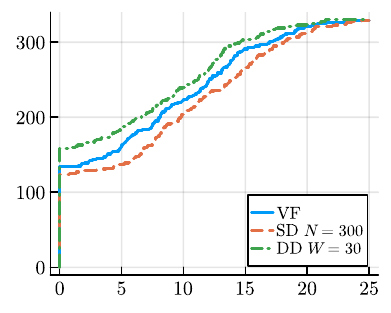}}
	}

	\subfloat[MaxSSPP - Time (s)]{
		\label{fig:exp-maxstab-time}
		\scalebox{0.85}{\includegraphics{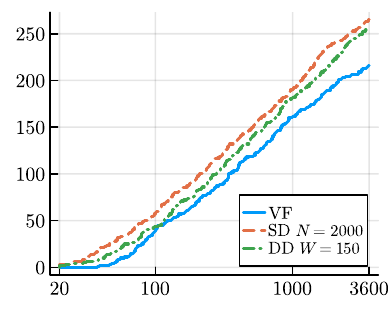}}
	}
	\subfloat[MinSCPP - Time (s)]{
		\label{fig:exp-mincover-time}
		\scalebox{0.85}{\includegraphics{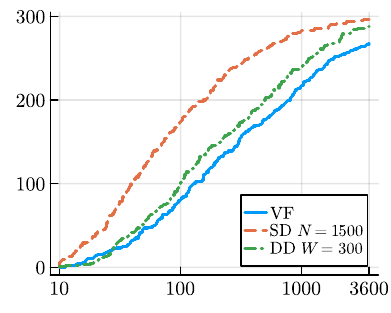}}
	}
 
    \subfloat[KIP - Time (s)]{
		\label{fig:exp-kip-time}
		\scalebox{0.85}{\includegraphics{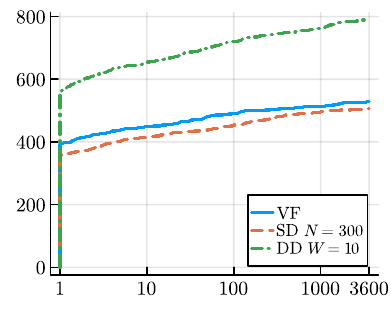}}
	}
	\subfloat[KIP - Optimality gap (\%)]{
		\label{fig:exp-kip-gap}
		\scalebox{0.85}{\includegraphics{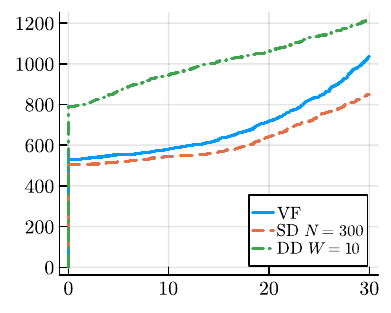}}
	}

	\caption{Cumulative number of instances ($y$-axis) with respect to solution time and optimality gap ($x$-axis).
		The time is depicted in logarithmic scale, while a linear scale is used for the gap.}
	\label{fig:exp-graphs}
\end{figure}

We tested the three DP models on the set of 1340 instances mentioned in \citep{weninger2023}.
In particular, the instances are taken from CCLW \citep{caprara2016}, DCS \citep{della2020},
DeNegre \citep{denegre2011}, FMS \citep{fischetti2018}, and TRS \citep{tang2016}.
However, we exclude the 1500 newly-generated instances in \citep{weninger2023} because they
have a very high number of items and our approaches are unable to solve them; we recall that our framework is aimed to be general-purpose.
The results are shown in \cref{tab:kip-results,fig:exp-kip-time,fig:exp-kip-gap}.
The decision diagram is clearly superior than the other two models in every class of instances.
Comparing the DP models to the state of the art method described in \citep{weninger2023},
our models are worse, but our experiments are executed on a single thread,
while theirs were run using 16 threads. Nonetheless, the DP models are more general and
should not be expected to outperform problem-specific algorithms. Analogously, the branch-and-bound method by~\cite{fischetti2018} for interdiction problems with monotonicity (recall~\cref{rm:kpp-maximal}) is overall faster than our models, for the instances they tested: CCLW, TRS and part of DeNegre. Despite this, our decision diagram model with $W=10$ solved all those instances in much less time than the time limit and it was significantly faster than the running times reported for the general bilevel solver MIX++ in~\cite{fischetti2017}. Although the main focus of this work are combinatorial pricing problems, these results allow us to provide further insights on the performance of our methodology.

\begin{table}[htbp]
	\centering
    \caption{Numerical results for the KIP.}
    \label{tab:kip-results}
    \vskip 2mm
	\begin{talltblr}[
        label = none,
		note{a} = {Geometric average. Values less than 1 are lifted to 1.},
		note{b} = {Arithmetic average.},
		note{} = {The best configuration of each model is highlighted in bold.},
		]{
		colspec = {lrrrrrrrrr},
		cells = {font=\small},
		cell{1}{-} = {c},
		cell{1}{1,2,3,Z} = {r=2}{},
		cell{1}{4} = {c=6}{},
		row{3,6,14} = {font=\bf\small},
			}
		\toprule
		{DP                                      \\ Model} & $N$/$W$ & {Num. of \\ instances \\ solved} & Time\TblrNote{a}~~(s) & & & & & & {Gap\TblrNote{b} \\ (\%)} \\
		\cmidrule[lr]{4-9}
		   &      &     & Total & CCLW & DCS & DeNegre & FMS & TRS &     \\
		\midrule
         VF &     0 &   528 &  197 &   10 & 2415 &    3 &  574 &    1 &    16 \\
        		\midrule
         SD &   100 &   497 &  236 &   18 & 2477 &    5 &  785 &    1 &    18 \\
         SD &   200 &   500 &  238 &   19 & 2421 &    5 &  799 &    1 &    18 \\
         SD &   300 &   506 &  239 &   22 & 2390 &    5 &  789 &    1 &    18 \\
         SD &   500 &   504 &  244 &   22 & 2454 &    6 &  788 &    1 &    18 \\
         SD &   700 &   506 &  245 &   27 & 2444 &    6 &  784 &    1 &    18 \\
         SD &  1000 &   501 &  252 &   29 & 2519 &    6 &  795 &    1 &    18 \\
         SD &  1500 &   503 &  252 &   27 & 2556 &    6 &  794 &    1 &    18 \\
         SD &  2000 &   502 &  255 &   28 & 2595 &    6 &  797 &    1 &    18 \\
         SD &  2500 &   500 &  257 &   30 & 2643 &    6 &  790 &    1 &    18 \\
         SD &  3000 &   503 &  258 &   29 & 2641 &    6 &  803 &    1 &    18 \\
        		\midrule
         DD &    10 &   790 &   51 &    2 &  346 &    1 &  154 &    1 &     8 \\
         DD &    20 &   786 &   54 &    2 &  378 &    1 &  168 &    1 &     8 \\
         DD &    30 &   783 &   56 &    3 &  398 &    1 &  170 &    1 &     8 \\
         DD &    50 &   780 &   59 &    3 &  423 &    1 &  177 &    1 &     9 \\
         DD &    70 &   770 &   62 &    3 &  443 &    1 &  184 &    1 &     9 \\
         DD &   100 &   763 &   64 &    3 &  458 &    2 &  187 &    1 &     9 \\
         DD &   150 &   759 &   66 &    4 &  469 &    2 &  187 &    1 &     9 \\
         DD &   200 &   756 &   68 &    4 &  482 &    2 &  184 &    1 &     9 \\
         DD &   250 &   751 &   69 &    5 &  483 &    2 &  184 &    1 &     9 \\
         DD &   300 &   741 &   70 &    5 &  477 &    2 &  181 &    1 &     9 \\
		\bottomrule
	\end{talltblr}
\end{table}

\subsection{Discussions}
\label{ssec:discussions}

The experiments on four different problems give us a well-rounded view on the effectiveness
of the diagram models. The decision diagram is robust and can provide a small boost in any circumstance
compared to the value function reformulation.
The selection diagram, however, only performs well in the MaxSSPP and the MinSCPP.
The central mechanism of the improvement is due to the reductions in the number of callback calls and the time spent per call.
In particular, most of the solution time is spent to solve the subproblem in the MaxSSPP and MinSCPP,
while for the KPP and the KIP, that time is spent mainly on branch-and-bound.

We remark the dissimilarities in the structures of the four problems.
In the KPP and the KIP, any pair of items can appear in a feasible solution (as long as $C$ is large enough).
This is not the case for the MaxSSPP, whose main constraint is the exclusion of certain pairs of items.
In the MinSCPP, although any pair can also appear in a feasible solution, the minimality of a cover discourages
the selection of pair of items that have many overlaps in $\setE$.
Furthermore, a maximal solution of the KPP contains around 50\% of all items.
This number is smaller in the other two problems: 20\% for the MaxSSPP and 10\% for the MinSCPP.
The disparity in the distribution of pairs in feasible solutions may explain
the unevenness in the performance of the selection diagram across the three problems. This analysis of item pairs for a given CPP appears to be a good indicator of whether~\cref{alg:sdiag-dynamic} used to solve the selection diagram model will be more effective than~\cref{alg:ddiag-dynamic} used to solve the decision diagram model. Recall that~\cref{alg:sdiag-dynamic} limits the number of layers in the selection diagram and the number of nodes in layer two. Thus, this analysis provides an empirical way to anticipate the potential strength of the (relaxed\footnotemark) dual representation in the selection diagram model and thus, of~\cref{alg:sdiag-dynamic}'s effectiveness.

\footnotetext{Recall that~\cref{alg:sdiag-dynamic} dynamically generates~\cref{cons:dp-single-arc}.}

The problems' structures may also explain why the KIP and KPP are effectively solved using selection diagram and decision diagram models with small values of $N$ and $W$. By contrast, these models achieve the best performances for larger values of $N$ and $W$ when used to solve the MaxSSPP and MinSCPP. A deeper understanding of the impact of $N$ and $W$ on our methodology's performance will require a theoretical study, particularly on the strength of the constraints generated through our solving procedure.

\section{Conclusion}
\label{sec:conclusion}

In this paper, we derived a single-level reformulation of the CPP
by taking the dual of a linear program corresponding to a dynamic programming model of the follower's problem,
and combining it with the original primal problem using strong duality.
We investigated two dynamic programming models: selection diagram and decision diagram,
and provided the algorithms to solve them dynamically.
Then, we tested these models in three specializations of the CPP as well as the KIP. We observed that the decision diagram model can consistently lead to performance improvements over the value function reformulation, while the selection diagram formulation only brings (significant) speedups for two  specializations of the CPP.

Future research directions include the exploration of other dynamic programming models
and the application of this technique to other problems in bilevel optimization.


\appendix

\section{Bounds for the McCormick envelope}
\label{app:bounds}

To solve \cref{prog:dp-single} as a MILP, one needs to linearize the bilinear terms
$\tx$ using the McCormick envelope (\cref{eq:mccormick}).
As explained in \cref{ssec:vf-form}, $M_i = v_i$ is a valid bound for the KPP
due to \cref{rm:kpp-t-bound}.
By a similar reasoning, $M_i = v_i$ is still a valid bound for the MaxSSPP.
However, this logic crumbles in the case of the MinSCPP because of the optimization direction.

To derive a valid bound for the MinSCPP, we take inspiration from the same linearization problem
occurring in the NPP \citep{dewez2008}.
Recall that in the MinSCPP, each $i \in \setI$ is a subset of $\setE$ (see \cref{sssec:dd-mincover}).
For $i \in \setI$, a valid bound $M_i$ for $t_ix_i$ is:
\[ \max \{0, \min \{ P_i, Q_i \} \} \]
where
\begin{align*}
    P_i &= f(\infty; i) - v_i, \\
    Q_i &= f(\infty; \setE) - f(0; \setE \setminus i) - v_i.
\end{align*}

The function $f(t; E)$ is defined to be
the total cost to cover all elements in $E \subseteq \setE$ given tolls $t$:
\[ f(t; E) = \min_x \{ \tr{(v + t)}x \mid A_Ex \geq \vecone,~x \in \setbin^{\setI} \}\]
where $A_E \in \{0, 1\}^{E \times \setI}$ is the incidence matrix between $E$ and $\setI$.
In particular, $f(\infty; E)$ is the cost to cover $E$ without using any tolled set in $\setI_1$
(also called the \emph{toll-free cost}),
while $f(0; E)$ is the cost to cover $E$ when all tolled sets are available for no extra costs
(also called the \emph{null-toll cost}).

The bound $P_i$ is valid because if $t_i > P_i$, then we can replace the set $i$
with the toll-free sets in the optimal solution of $f(\infty; i) < v_i + t_i$,
which still cover all elements in $i$ but with less total cost.

The bound $Q_i$ is valid because the lowest cost for any solution that includes $i$
is $f(0; \setE \setminus i) + v_i + t_i$. If $t_i > Q_i$,
then the toll-free solution of $f(\infty; \setE)$ will have lower cost than any solution
that contains $i$, meaning the latter can never be optimal when $t_i > Q_i$.

\section{Effectiveness of Item Grouping}
\label{app:eff-grouping}

To prove that the item grouping technique described in \cref{ssec:grouping} is beneficial to the decision diagram models, we ran additional experiments comparing models with and without item grouping. The values of the width $W$ were chosen to match those of the best decision diagram models in \cref{sec:exp}. For MaxSSPP and MinSCPP, we also tested intermediate $W$ values due to their similar performance. The results are shown in \cref{tab:grouping-results}, which confirms that item grouping improves the performance of decision diagram models in general, with the exception of the KIP. We note that the best decision diagram model for the KIP in \cref{tab:kip-results} has $W = 10$. This suggests that in the case of the KIP, it is preferable to have more depth (no item grouping) rather than more width (higher $W$).

 These results are based on the same random seeds (recall \cref{alg:dd-initial}), reducing variance in relative performance. However, the seeds are different than those used in \cref{sec:exp}, likely leading to distinct initial decision diagrams, which explains the slight performance differences (with item grouping).

\begin{table}[htbp]
	\centering
    \caption{Performance of DD models with and without item grouping.}
    \label{tab:grouping-results}
    \vskip 2mm
	\begin{talltblr}[
        label = none,
		note{a} = {Geometric average.},
		note{b} = {Arithmetic average.},
		note{} = {The best configuration of each pair is highlighted in bold.},
		]{
		colspec = {lrrrrr},
		cells = {font=\small},
		cell{1}{-} = {c,m},
        cell{2}{1} = {r=2}{},
        cell{4}{1} = {r=6}{},
        cell{10}{1} = {r=4}{},
        cell{14}{1} = {r=2}{},
		cell{2,4,6,8,10,12,15}{2-Z} = {font=\bf\small},
		}
		\toprule
		{Problem} & {Item \\ grouping} & $W$ & {Num. of \\ instances \\ solved} & Time\TblrNote{a}~~(s) & {Gap\TblrNote{b} \\ (\%)} \\
		\midrule
        KPP     & Yes &    30 &   161 &    986 &   5.2 \\
                &  No &    30 &   147 &   1160 &   5.6 \\
        \midrule
        MaxSSPP & Yes &    10 &   252 &    484 &   4.9 \\
                &  No &    10 &   250 &    494 &   4.5 \\
        \cmidrule[lr]{2-Z}
                & Yes &    50 &   259 &    413 &   4.3 \\
                &  No &    50 &   251 &    535 &   5.1 \\
        \cmidrule[lr]{2-Z}
                & Yes &   150 &   258 &    444 &   2.6 \\
                &  No &   150 &   242 &    737 &   7.8 \\
        \midrule
        MinSCPP & Yes &    50 &   281 &    238 &   5.6 \\
                &  No &    50 &   279 &    252 &   7.8 \\
        \cmidrule[lr]{2-Z}
                & Yes &   300 &   284 &    214 &   4.0 \\
                &  No &   300 &   250 &    911 &  15.9 \\
        \midrule
        KIP     & Yes &    10 &   800 &   49.6 &   8.0 \\
                &  No &    10 &   880 &   29.2 &   7.8 \\
		\bottomrule
	\end{talltblr}
\end{table}

\section{Sensitivity Analysis of the Instance Generation Algorithms}
\label{app:sens-anal}
During our preliminary analysis on the randomly generated instances, we discovered that several combinations of generation parameters yielded easier-to-solve instances. Since our aim is to improve the time required to solve the hard instances, we wish to tune the generation algorithm to avoid producing trivial instances. Thus, we optimized the difficulty of the generated instances with the CMA-ES algorithm \citep{hansen2003}. However, the difficulty is not evaluated directly (\emph{i.e.}, by solving the instances and measuring the solution time) since it takes too much time. Instead, it is estimated as $[f(0) - f(\infty)]/g$, where $f(t)$ is the maximum follower's objective value given $t$:
\begin{equation*}
	\label{eq:follower-val}
	f(t) = \max_x \{ \vtx \mid x \in \setX \}
\end{equation*}
and $g$ is the maximum leader's revenue:
\begin{equation*}
	g = \max_{t,x} \{ \tx \mid t \in \setT,~x \in \funcR(t) \}.
\end{equation*}

The term $f(\infty)$ represents the follower's utility without the tolled items while $f(0)$ represents the highest utility that the follower can ever obtain. The leader's revenue is bounded by the extra utility that the leader's items can provide, hence $f(0) - f(\infty) \geq g$. We speculate that a CPP instance is hard when the bound $f(0) - f(\infty)$ on the revenue is high, but the actual revenue $g$ is low, because, in this case, the follower has more complex reactions rather than two trivial reactions: one for $t_i = 0$ for all $i \in \setIa$, and one for $t_i \to \infty$ for all $i \in \setIa$. Using this estimate, we can optimize the generation parameters using the CMA-ES algorithm, while generating instances with a lower number of items $N$ for faster evaluation.

\begin{figure}[htbp]
	\centering
	\scalebox{0.9}{\includegraphics{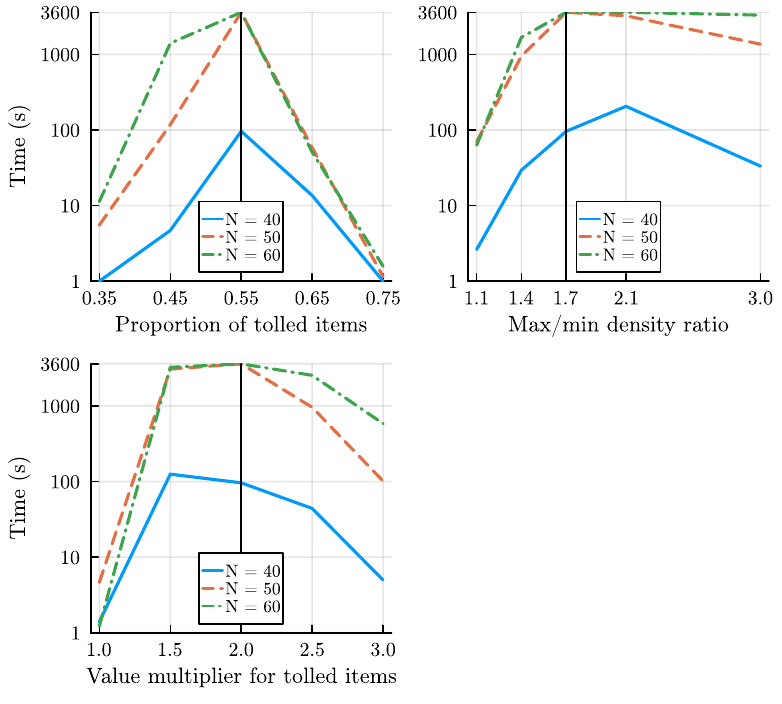}}
	\caption{Geometric average time with respect to various parameter values for the KPP instance generation. Longer time means harder instances. The thick vertical lines mark the baseline values.}
	\label{fig:kpp-params}
\end{figure}

\begin{figure}[htbp]
	\centering
	\scalebox{0.9}{\includegraphics{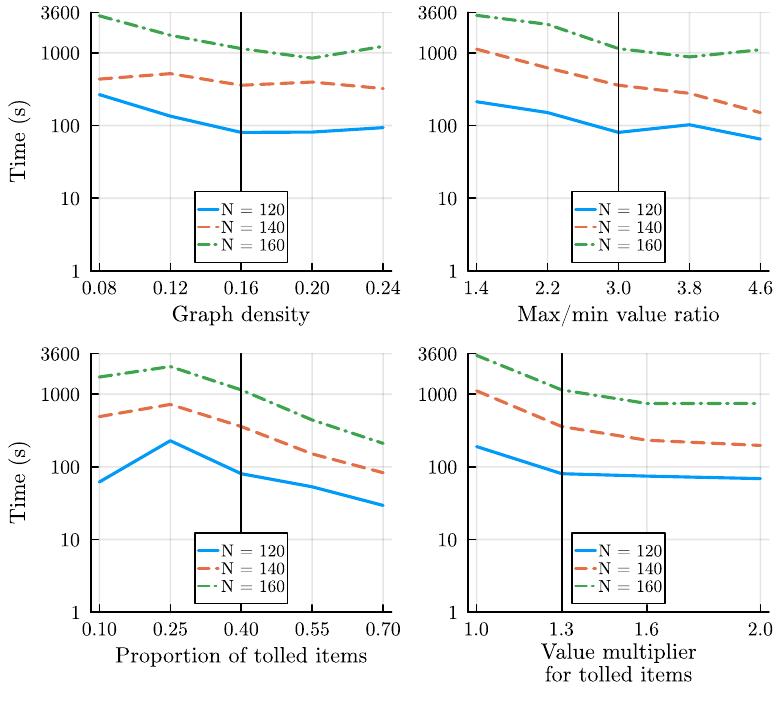}}
	\caption{Geometric average time with respect to various parameter values for the MaxSSPP instance generation. Longer time means harder instances. The thick vertical lines mark the baseline values.}
	\label{fig:maxstab-params}
\end{figure}

\begin{figure}[htbp]
	\centering
	\scalebox{0.9}{\includegraphics{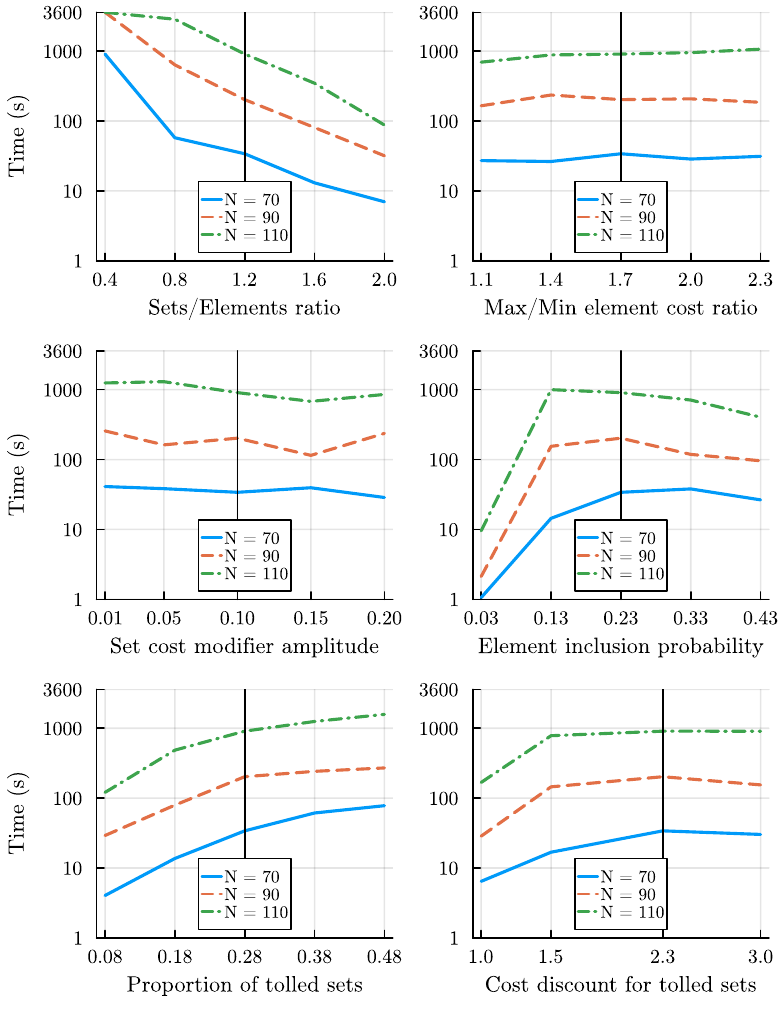}}
	\caption{Geometric average time with respect to various parameter values for the MinSCPP instance generation. Longer time means harder instances. The thick vertical lines mark the baseline values.}
	\label{fig:mincover-params}
\end{figure}

In this appendix, we want to verify that the set of parameters used to generate the instances for the main experiments did not produce overly-trivial instances. We varied a single parameter from a baseline set of parameters (which is the set used in the main experiments), then generated 10 instances, measured the average solution time of the value-function model, and compared it to the baseline result. The results are shown in \cref{fig:kpp-params,fig:maxstab-params,fig:mincover-params}.

The results show that the baseline set of parameters, although it does not yield the hardest instances, also does not produce trivial instances. The KPP instances are specially sensitive to the parameter values, with the solution time dramatically dropping at the two ends of the graphs. The MaxSSPP and MinSCPP are more resilient to the generation parameters. Overall, we conclude that the baseline set of parameters is adequate for the main experiments, although there is more room to increase the difficulty of the instances. These results also give insights into what make a particular CPP hard, especially when CPPs are not yet well-researched.

\section*{Acknowledgments}
This work was funded by FRQ-IVADO Research Chair in Data Science for Combinatorial Game Theory, and the NSERC grant 2019-04557 and 2024-04051. 

This research was enabled in part by support provided by Calcul Qu\'ebec (\url{https://www.calculquebec.ca})
and the Digital Research Alliance of Canada (\url{https://alliancecan.ca}).

\bibliographystyle{plainnat}
\bibliography{ref}

\end{document}